\documentclass[11pt]{article}
\usepackage{amsmath}
\usepackage{amsfonts}
\usepackage{amssymb}
\usepackage{amsthm}
\usepackage{amscd}
\usepackage[vcentermath]{youngtab}

\newtheorem{theorem}{Theorem}

\newtheorem{corollary}[theorem]{Corollary}

\newtheorem{lemma}[theorem]{Lemma}

\newtheorem{proposition}[theorem]{Proposition}

\theoremstyle{definition}
\newtheorem{remark}[theorem]{Remark}
\newtheorem{example}[theorem]{Example}

\newcommand{\ba}{\overline{1}}
\newcommand{\bbb}{\overline{2}}
\newcommand{\bc}{\overline{3}}
\newcommand{\bd}{\overline{4}}
\newcommand{\bn}{\overline{n}}
\newcommand{\bnm}{\overline{n\!\!-\!\!1}}

\newcommand{\Bh}{\widehat{B}}

\newcommand{\bor}{\mathfrak{b}}

\newcommand{\CC}{\mathcal{C}}

\newcommand{\Db}{\overline{D}}

\newcommand{\Dem}{E}
\newcommand{\DW}{\Z_{\ge}^n}
\newcommand{\dia}{\diamondsuit}
\newcommand{\et}{\widetilde{e}}
\newcommand{\gf}{\mathfrak{g}}
\newcommand{\F}{F}
\newcommand{\fb}[1]{\fbox{$#1$}}
\newcommand{\fbs}[1]{\fbox{\scriptsize{$#1$}}}

\newcommand{\ft}{\widetilde{f}}
\newcommand{\gaf}{\widehat{\mathfrak{g}}}

\newcommand{\Hb}{\overline{H}}

\newcommand{\hh}{\mathfrak{h}}
\newcommand{\I}{I}
\newcommand{\inner}[2]{\langle#1\,,\,#2\rangle}
\newcommand{\J}{J}
\newcommand{\Kb}{\overline{K}}
\newcommand{\KK}{\widetilde{D}}
\newcommand{\KL}{\mathrm{KL}}
\newcommand{\KLS}{{}^\infty\mathrm{KL}}
\newcommand{\la}{\lambda}

\newcommand{\lah}{\widehat{\la}}
\newcommand{\lwt}{L}
\newcommand{\Modf}{\mathrm{Mod}}
\newcommand{\muh}{\widehat{\mu}}

\newcommand{\Par}{\mathcal{P}}

\newcommand{\TT}{\mathbf{T}}
\newcommand{\tW}{\widetilde{W}}
\newcommand{\ve}{\varepsilon}
\newcommand{\vn}{\varnothing}
\newcommand{\vp}{\varphi}
\newcommand{\wt}{\mathrm{wt}}
\newcommand{\Xb}{\overline{X}}
\newcommand{\Z}{\mathbb{Z}}

\begin{document}

\title{Lusztig's $q$-analogue of weight multiplicity
and one-dimensional sums for affine root systems}

\author{C\'edric Lecouvey \\
Laboratoire de Mathématiques Pures et Appliqu\'ees Joseph Liouville \\
Universit\'e du Littoral \\
Centre Universitaire de la Mi-Voix \\
Maison de la Recherche Blaise Pascal \\
50 rue F.Buisson B.P. 699 62228 Calais Cedex \\ \\
Mark Shimozono\footnote{MS was supported in part by NSF grant DMS-0401012.} \\
Department of Mathematics \\
Virginia Polytechnic Institute and State University \\
Blacksburg, VA 24061-0123 USA }



\maketitle
\begin{abstract}
In this paper we complete the proof of the $X=K$ conjecture, that
for every family of nonexceptional affine algebras, the graded
multiplicities of tensor products of ``symmetric power''
Kirillov-Reshetikhin modules known as one-dimensional sums, have a
large rank stable limit $X$ that has a simple expression (called the
$K$-polynomial) as nonnegative integer combination of Kostka-Foulkes
polynomials. We consider a subfamily of Lusztig's $q$-analogues of
weight multiplicity which we call stable KL polynomials and denote
by $\KLS$. We give a type-independent proof that $K=\KLS$. This
proves that $X=\KLS$: the family of stable one-dimensional sums
coincides with family of stable KL polynomials. Our result
generalizes the theorem of Nakayashiki and Yamada which establishes
the above equality in the case of one-dimensional sums of affine
type A and the Lusztig $q$-analogue of type $A$, where both are
Kostka-Foulkes polynomials.
\end{abstract}


\section{Introduction}

One-dimensional (1-d) sums $X$ are graded tensor product
multiplicities for affine Kac-Moody algebras, which arise from
two-dimensional solvable lattice models \cite{KMN} and which may be
defined using the combinatorics of affine crystal graphs
\cite{HKOTY} \cite{HKOTT}. The definition of a 1-d sum depends not
only on an affine algebra $\gaf$ but also on a distinguished simple
Lie subalgebra $\gf$. This is the same data as giving an affine
Dynkin diagram together with a distinguished $0$ node whose removal
leaves the Dynkin diagram of a simple Lie algebra.

Within any nonexceptional family of affine algebras, the 1-d sums
have a large rank limit that we shall call stable 1-d sums. There
are many nonexceptional families of affine algebras and choices of
$0$ node, but only four distinct kinds of stable 1-d sums \cite{SZ}
\cite{sh}. They depend not on the entire affine Dynkin diagram, but
only on the neighborhood of the $0$ node. The four kinds of stable
1-d sums $X^\dia$ are labeled by the four partitions $\dia\in\{\vn,
(1), (2), (1,1)\}$ having at most 2 cells. This labeling is inspired
by the branching rules for the restriction of Kirillov-Reshetikhin
modules \cite{HKOTY} \cite{HKOTT} from the affine algebra to its
simple Lie subalgebra.

In \cite{SZ} it was conjectured that every stable 1-d sum $X$ has a
surprisingly simple explicit expression (called the $K$ polynomial)
in terms of those of type $A$, the latter being the well-known
Kostka-Foulkes polynomials. We finish the proof of this $X=K$
conjecture for tensor products of symmetric powers, supplying the
proof in the case $\dia=(1,1)$; the other cases were already settled
in \cite{sh}. Unlike the situation in \cite{sh}, for the case
$\dia=(1,1)$ we cannot use affine crystal theory alone, but must use
a novel grading that is not compatible with the affine crystal
structure but only with an $A_{n-1}$ crystal substructure.

On the other hand there is another family of polynomials, Lusztig's
$q$-analogue of weight multiplicity for classical Lie algebras
\cite{Lu}. These are instances of affine Kazhdan-Lusztig
polynomials. There is a subfamily which we shall call stable KL
polynomials and denote by $\KLS$, which are precisely those Lusztig
$q$-analogues that are indexed by pairs of dominant weights that
contain sufficiently large multiples of a certain fundamental
weight. They are called stable because if both weights are
positively translated by the above fundamental weight, then the
polynomial does not change. The stable KL polynomials were studied
in \cite{lec} \cite{lec4}.

We give a type-independent proof that $K=\KLS$. The stable 1-d sum
types $\vn,(1),(2),(1,1)$ correspond to stable KL polynomials of
types $A_{n-1},B_n,C_n,D_n$ respectively. For type $A$, $K=\KLS$
holds by definition. For types $C$ and $D$ the proof was given in
\cite{lec} \cite{lec4}. It was also noticed in \cite{lec4} that the
usual type $B$ Lusztig $q$-analogue does not agree with $X^{(1)}$,
despite the fact that both polynomials have the same value at $q=1$.
Here we repair this problem by using a $q$-analogue of Kostant's
partition function that uses different powers of $q$ for different
root lengths, which is related to Hecke algebras with unequal
parameters.

Combined with the $X=K$ theorem we obtain $X=\KLS$: the stable 1-d
sums are the stable KL polynomials.

In the $X=\KLS$ theorem the correspondence of weights involves a
twist $\mu\mapsto\muh$ that also appears in Howe duality \cite{H}.
We do not have an adequate explanation for this twist. Our proof
reduces to a type $A_{n-1}$ situation where the twist is given by
contragredient duality, but its meaning for the original weights of
classical type is less clear to us.

\subsection{Classical algebras and crystal graphs} Let $\gf\supset
\bor\supset\hh$ be a simple Lie algebra, a Borel subalgebra and a
Cartan subalgebra, $\J$ the set of nodes for the Dynkin diagram of
$\gf$, $\{\alpha_i\mid i\in\J\}\subset\hh^*$ the simple roots,
$R^+\subset\hh^*$ the set of positive roots with respect to $\bor$,
$\{h_i\mid i\in\J\}\subset\hh$ the simple coroots, $\{\omega_i\mid
i\in \J\}\subset\hh^*$ the fundamental weights, $P=\bigoplus_{i\in
\J} \Z \omega_i\subset\hh^*$ the weight lattice, $P^+ =
\bigoplus_{i\in \J} \Z_{\ge0} \omega_i$ the set of dominant weights,
$W$ the Weyl group, and $\inner{\cdot}{\cdot}$ the evaluation
pairing of $\hh$ with $\hh^*$. Let $U_q(\gf)$ be the quantized
universal enveloping algebra of $\gf$ \cite{Ka}.

Each finite-dimensional $U_q(\gf)$-module has a crystal graph
\cite{Ka}. Let $\CC(\gf)$ be the category of such crystal graphs.
Each $B\in\CC(\gf)$ is a directed graph with vertex set also denoted
$B$ and directed edges labeled by the set $\J$. If we remove all
edges of $B$ except those labeled by a fixed $i\in\J$, then the
connected components are directed paths called $i$-strings. Denote
by $\ft_i(b)$ (resp. $\et_i(b)$) be the next (resp. previous) vertex
in the $i$-string through the vertex $b\in B$. If no such vertex
exists then the result is defined to be the special symbol
$\emptyset$. Define $\et_i(\emptyset)=\ft_i(\emptyset)=\emptyset$.
Let $\ve_i(b)$ (resp. $\vp_i(b)$) be the number of steps to the
beginning (resp. end) of the $i$-string through $b$. Each
$B\in\CC(\gf)$ also has a weight function $\wt:B\rightarrow P$. It
satisfies
\begin{alignat*}{2}
\wt(\ft_i(b))&=\wt(b)-\alpha_i &\qquad&\text{if
$\ft_i(b)\not=\emptyset$} \\
\wt(\et_i(b))&=\wt(b)+\alpha_i &\qquad&\text{if
$\et_i(b)\not=\emptyset$} \\
\inner{h_i}{\wt(b)} &= \vp_i(b) - \ve_i(b) &&
\end{alignat*}
A vertex $b\in B$ such that $\ve_i(b)=0$ for all $i\in\J$, is called
a highest weight vertex of $B$. For $B,B'\in\CC(\gf)$ a morphism
$\psi:B\rightarrow B'$ is a map $\psi:B\cup\{\emptyset\}\rightarrow
B'\cup\{\emptyset\}$ such that $\psi(\emptyset)=\emptyset$ and
$\psi(\et_i(b))=\et_i(\psi(b))$ and $\psi(\ft_i(b))=\ft_i(\psi(b))$
for all $b\in B$ and $i\in\J$. An isomorphism is a bijective
morphism whose inverse function is also a morphism.

If $V$ and $V'$ are finite-dimensional $U_q(\gf)$-modules with
crystal graphs $B$ and $B'$, then their direct sum $V\oplus V'$ has
crystal graph $B \oplus B'$, which by definition is the disjoint
union of $B$ and $B'$.

Each $B\in\CC(\gf)$ is the disjoint union of connected crystal
graphs. For every connected $B\in \CC(\gf)$, there is a unique
$\la\in P^+$ such that $B\cong B(\la)$ where $B(\la)$ is the crystal
graph of the irreducible finite-dimensional $U_q(\gf)$-module
$V(\la)$ of highest weight $\la$. For $\gf$ of classical type
$B(\la)$ is described explicitly in \cite{KN}. $B(\la)$ has a unique
highest weight vertex, which is also the unique vertex of weight
$\la$ in $B(\la)$.

The tensor product $V\otimes V'$ has the crystal graph $B \otimes
B'$ whose underlying set is the Cartesian product $B \times B'$,
whose elements are written $b\otimes b'$ instead of $(b,b')$, with
crystal structure given by
\begin{equation} \label{E:tensor}
\ft_i(b\otimes b') = \begin{cases}%
  \ft_i(b) \otimes b' & \text{if $\vp_i(b)>\ve_i(b')$} \\ %
   b\otimes \ft_i(b') & \text{if $\vp_i(b)\le \ve_i(b')$}
\end{cases} \qquad
\et_i(b\otimes b') = \begin{cases}%
  \et_i(b) \otimes b' & \text{if $\vp_i(b) \ge \ve_i(b')$} \\ %
   b\otimes \et_i(b') & \text{if $\vp_i(b) < \ve_i(b')$.}
\end{cases}
\end{equation}

The following lemma is a straightforward consequence of
\eqref{E:tensor}.

\begin{lemma}
\label{lem_plu_hp} Let $B,B'\in\CC(\gf)$. Then $b\otimes b'\in
B\otimes B'$ is a highest weight vertex if and only if $b\in B$ is a
highest weight vertex and $\vp_i(b)\ge\ve_i(b')$ for any $i\in \J$.
\end{lemma}

If $b\in B'$ is a highest weight vertex in $B'\in\CC(\gf)$, we
denote by $B(b)$ the connected component of $b$ in $B'$.

\subsection{Partitions and dominant weights}
{}From now on let $\gf=\gf_n$ be a classical Lie algebra, one of
type $A_{n-1}$, $B_n$, $C_n$, or $D_n$, with Dynkin diagrams labeled
as in Figure \ref{fig:rootsys}. We shall also use a nonstandard
labeling $D_n^\dagger$.
\begin{figure}
\begin{equation*}
\begin{tabular}[t]{|c|l|} \hline
\begin{minipage}[b]{4em}
\begin{center}
$A_{n-1}$ \\ $(n \ge 2)$  \\[1mm]
\end{center}
\end{minipage}&
\begin{picture}(106,30)(-5,-15)
\multiput( 0,0)(20,0){2}{\circle{6}}%
\multiput(80,0)(20,0){2}{\circle{6}}%
\multiput( 3,0)(20,0){2}{\line(1,0){14}}%
\multiput(63,0)(20,0){2}{\line(1,0){14}}%
\multiput(39,0)(4,0){6}{\line(1,0){2}}%
\put(0,-5){\makebox(0,0)[t]{\scriptsize{$1$}}} %
\put(20,-5){\makebox(0,0)[t]{\scriptsize{$2$}}}%
\put(80,-5){\makebox(0,0)[t]{\scriptsize{$n\!\! -\!\! 2$}}}%
\put(100,-5){\makebox(0,0)[t]{\scriptsize{$n\!\! -\!\! 1$}}}%
\end{picture}
\\ \hline
\begin{minipage}[b]{4em}
\begin{center}
$B_n$\\$(n \ge 2)$ \\[1mm]
\end{center}
\end{minipage}&
\begin{picture}(130,30)(-5,-15)
\multiput( 0,0)(20,0){2}{\circle{6}}%
\multiput(80,0)(20,0){3}{\circle{6}}%
\multiput( 3,0)(20,0){2}{\line(1,0){14}}%
\multiput(63,0)(20,0){2}{\line(1,0){14}}%
\multiput(39,0)(4,0){6}{\line(1,0){2}}%
\put(0,-5){\makebox(0,0)[t]{\scriptsize{$1$}}} %
\put(20,-5){\makebox(0,0)[t]{\scriptsize{$2$}}}%
\put(80,-5){\makebox(0,0)[t]{\scriptsize{$n\!\! -\!\! 2$}}}%
\put(100,-5){\makebox(0,0)[t]{\scriptsize{$n\!\! -\!\! 1$}}}%
\put(120,-6){\makebox(0,0)[t]{\scriptsize{$n$}}}%
\multiput(102.85,-1)(0,2){2}{\line(1,0){14.3}} 
\put(110,0){\makebox(0,0){$>$}}%
\end{picture}
\\ \hline
\begin{minipage}[b]{4em}
\begin{center}
$C_n$\\$(n \ge 2)$\\[1mm]
\end{center}
\end{minipage}&
\begin{picture}(130,30)(-5,-15)
\multiput( 0,0)(20,0){2}{\circle{6}}%
\multiput(80,0)(20,0){3}{\circle{6}}%
\multiput( 3,0)(20,0){2}{\line(1,0){14}}%
\multiput(63,0)(20,0){2}{\line(1,0){14}}%
\multiput(39,0)(4,0){6}{\line(1,0){2}}%
\put(0,-5){\makebox(0,0)[t]{\scriptsize{$1$}}} %
\put(20,-5){\makebox(0,0)[t]{\scriptsize{$2$}}}%
\put(80,-5){\makebox(0,0)[t]{\scriptsize{$n\!\! -\!\! 2$}}}%
\put(100,-5){\makebox(0,0)[t]{\scriptsize{$n\!\! -\!\! 1$}}}%
\put(120,-6){\makebox(0,0)[t]{\scriptsize{$n$}}}%
\multiput(102.85,-1)(0,2){2}{\line(1,0){14.3}} 
\put(110,0){\makebox(0,0){$<$}}%
\end{picture}
\\ \hline
\begin{minipage}[b]{4em}
\begin{center}
$D_n$ \\ $(n \ge 4)$ \\ {\hphantom{x}}
\end{center}
\end{minipage}&
\begin{picture}(110,45)(-5,-15)
\multiput( 0,0)(20,0){2}{\circle{6}}
\multiput(80,0)(20,0){2}{\circle{6}}%
\put(80,20){\circle{6}} %
\multiput(3,0)(20,0){2}{\line(1,0){14}} %
\multiput(63,0)(20,0){2}{\line(1,0){14}}
\multiput(39,0)(4,0){6}{\line(1,0){2}}%
\put(80,3){\line(0,1){14}}%
\put(0,-5){\makebox(0,0)[t]{\scriptsize{$1$}}} %
\put(20,-5){\makebox(0,0)[t]{\scriptsize{$2$}}}
\put(80,-5){\makebox(0,0)[t]{\scriptsize{$n\!\! - \!\! 2$}}}
\put(100,-5){\makebox(0,0)[t]{\scriptsize{$n\!\! -\!\! 1$}}}
\put(85,20){\makebox(0,0)[l]{\scriptsize{$n$}}}
\end{picture}
\\ \hline
\begin{minipage}[b]{4em}
\begin{center}
$D_n^\dagger$ \\ $(n \ge 4)$ \\ {\hphantom{x}}
\end{center}
\end{minipage}&
\begin{picture}(110,45)(-5,-15)
\multiput( 0,0)(20,0){2}{\circle{6}}
\multiput(80,0)(20,0){2}{\circle{6}}%
\put(80,20){\circle{6}} %
\multiput(3,0)(20,0){2}{\line(1,0){14}} %
\multiput(63,0)(20,0){2}{\line(1,0){14}}
\multiput(39,0)(4,0){6}{\line(1,0){2}}%
\put(80,3){\line(0,1){14}}%
\put(0,-5){\makebox(0,0)[t]{\scriptsize{$n\!\!-\!\!1$}}} %
\put(20,-5){\makebox(0,0)[t]{\scriptsize{$n\!\!-\!\!2$}}}
\put(80,-5){\makebox(0,0)[t]{\scriptsize{$2$}}}
\put(100,-5){\makebox(0,0)[t]{\scriptsize{$1$}}}
\put(85,20){\makebox(0,0)[l]{\scriptsize{$0$}}}
\end{picture}
\\ \hline
\end{tabular}
\end{equation*}
\caption{\label{fig:rootsys} Dynkin diagrams for classical Lie
algebras}
\end{figure}
Let $\J$ be the vertex set of the Dynkin diagram of $\gf_n$:
$\J=\{1,2,\dotsc,n-1\}$ for $\gf_n=A_{n-1}$ and
$\J=\{1,2,\dotsc,n\}$ for $\gf_n=B_n,C_n,D_n$. The diagram
$D_n^\dagger$ has vertex set $\J^\dagger=\{0,1,\dotsc,n-1\}$.

The weight lattice $P$ of $\gf_n$ may be explicitly realized as a
sublattice of $(\frac{1}{2}\Z)^n$ where $\omega_i=(1^i,0^{n-i})$ for
$1\le i\le n-2$ and also $i=n-1$ for $\gf_n\not=D_n$,
$\omega_n^{C_n}=(1^n)$,
$\omega_n^{B_n}=\omega_n^{D_n}=(\frac{1}{2}^n)$, and
$\omega_{n-1}^{D_n}=(\frac{1}{2}^{n-1},-\frac{1}{2})$. For
$D_n^\dagger$, $P$ is identified with the sublattice of
$(\frac{1}{2}\Z)^n$ generated by the fundamental weights
$\omega^{D_n^\dagger}_0 = ((-\frac{1}{2})^n)$,
$\omega^{D_n^\dagger}_1 = (\frac{1}{2},(-\frac{1}{2})^{n-1})$, and
$\omega^{D_n^\dagger}_i=(0^i,(-1)^{n-i})$ for $2\le i\le n-1$.

The sum of fundamental weights (or half-sum of positive roots)
$\rho$ may be given by $\rho^{A_{n-1}}=(n-1,n-2,\dotsc,1,0)$,
$\rho^{B_n}=\frac{1}{2}(2n-1,2n-3,\dotsc,1)$,
$\rho^{C_n}=(n,n-1,\dotsc,2,1)$, and
$\rho^{D_n}=(n-1,n-2,\dotsc,1,0)$.

Let $\ve_i$ be the $i$-th standard basis vector. The simple roots
$\{\alpha_i\mid i\in J \}\subset \hh^*$ are given by
$\alpha_i=\ve_i-\ve_{i+1}$ for $1\le i\le n-1$,
$\alpha_n^{B_n}=\ve_n$, $\alpha_n^{C_n}=2\ve_n$,
$\alpha_n^{D_n}=\ve_{n-1}+\ve_n$, and
$\alpha_0^{D_n^\dagger}=-\ve_1-\ve_2$.

Later we shall make use of the fact that the root system of type
$A_{n-1}$ is the subsystem of $C_n$ (resp. $D_n^\dagger$) that is
obtained by removing the Dynkin node $n$ (resp. $0$).

We identify elements of $\Z^n \subset P$ with weights (excluding the
case of $D_n^\dagger$ here). The set $P^+ \cap \Z^n$ of dominant
weights in the sublattice $\Z^n$ is given by partitions
$\la=(\la_1\ge\la_2\ge\dotsm\ge \la_n)\in\Z^n$ where $\la_n\ge0$
except in type $D_n$ in which case $\la_{n-1} \ge|\la_n|$. If $\la$
is a partition with at most $m$ parts and $n\ge m$ then we may view
$\la$ as a dominant weight for $\gf_n$.

\subsection{Affine algebras, KR crystal graphs, and 1-d sums}
Let $\gaf\supset \gaf'\supset \gf$ be an affine Kac-Moody algebra
\cite{Kac}, its derived subalgebra, and the simple Lie subalgebra
whose Dynkin diagram is obtained from the affine Dynkin diagram by
removing a distinguished zero node. Let $U_q(\gaf)\supset
U'_q(\gaf)\supset U_q(\gf)$ be the corresponding Drinfeld-Jimbo
quantized universal enveloping algebras \cite{Ka}.

Let $U'_q(\gaf)-\Modf$ be the category of finite-dimensional
irreducible $U'_q(\gaf)$-crystals which possess an affine crystal
graph\footnote{For finite-dimensional $U'_q(\gaf)$-modules this is
rare. In comparison, the nonzero irreducible integrable highest
weight $U'_q(\gaf)$-modules are all known to have affine crystal
graphs, but all these modules are infinite-dimensional.}. Such
affine crystal graphs have directed edges colored by the set
$\I=\J\cup\{0\}$ and weight function with values in $P$ (after
projection from the weight lattice of $\gaf'$). It was conjectured
in \cite{HKOTY} \cite{HKOTT} that for every $(r,s)\in \J \times
\Z_{>0}$ there is a module $W^{(r)}_s\in U'_q(\gaf)-\Modf$  called
the Kirillov-Reshetikhin (KR) module, with affine crystal graph
denoted $B^{r,s}$. The KR modules are fundamental objects: it is
expected that every $V\in U'_q(\gaf)-\Modf$ is isomorphic to a
tensor product of KR modules. The KR module $W_s^{(r)}$ has a
prescribed $U_q(\gf)$-decomposition of the form \cite{HKOTY}
\cite{HKOTT}
$$W^{(r)}_s \cong V(s\omega_r) \oplus \dotsm$$ where
$\dotsm$ indicates a direct sum of irreducibles $V(\mu)$ where $\mu$
is smaller than $s\omega_r$ in a certain sense.

{}From now on we assume that $\gaf$ is a nonexceptional affine
algebra. The ``symmetric power" KR modules $W_s^{(1)}$ and their
affine crystal graphs $B^{1,s}$ were constructed in \cite{KMN2}
\cite{KKM}. Let $\CC$ be the category of tensor products of KR
crystals of the form $B^{1,s}$. Let
\begin{equation} \label{E:d}
  d = \begin{cases}
  2 & \text{if the Dynkin diagram of $\gaf$ has a double bond
  directed from node $1$ to node $0$} \\
  1 & \text{otherwise.}
  \end{cases}
\end{equation}
$\CC$ has the following remarkable properties.
\begin{enumerate}
\item Every $B\in\CC$ is the affine crystal graph of an irreducible
$U'_q(\gaf)$-module. In particular $B$ is connected as an affine
crystal graph.
\item For every $B_1,B_2\in\CC$ there is a unique affine crystal
isomorphism $R_{B_1,B_2}:B_1 \otimes B_2 \rightarrow B_2 \otimes
B_1$ called the combinatorial $R$-matrix, and a map
$\Hb=\Hb_{B_1,B_2}:B_1\otimes B_2\rightarrow \frac{1}{d}\Z$ called
the local coenergy function, such that the following holds.
$b_1\otimes b_2\in B_1 \otimes B_2$ and $b_2'\otimes
b_1'=R_{B_1,B_2}(b_1\otimes b_2)$. Then
\begin{gather}
\text{For $b_1\otimes b_2\in B_1\otimes B_2$ and $b_2'\otimes
b_1'=R_{B_1,B_2}(b_1\otimes b_2)$,} \hphantom{xxxxxxxxxxxxxxxxxxxx} \nonumber \\
\label{E:H} %
\Hb(\et_i(b_1\otimes b_2))=\Hb(b_1\otimes
b_2)+\dfrac{\delta_{i0}}{d}
\begin{cases}
1 & \text{if $\vp_{0}(b_1)\ge\ve_{0}(b_2)$ and
$\vp_{0}(b_2')\ge\ve_{0}(b_1')$,} \\ %
-1 & \text{if $\vp_{0}(b_1)<\ve_{0}(b_2)$ and
$\vp_{0}(b_2')<\ve_{0}(b_1')$,} \\ %
0 & \text{otherwise.}
\end{cases}
\end{gather}
Thus $\Hb$ is a $U_q(\gf)$-equivariant grading of $B_1\otimes B_2$
(that is, it is constant under directed edges in $\J$). Since
$B_1\otimes B_2\in\CC$ is connected as an affine crystal graph, the
map $\Hb$ is unique up to a global additive constant. If
$B_1=B_2=B$, we have $R_{B,B}=1_{B\otimes B}$ by uniqueness and
\begin{equation} \label{E:H=}
\Hb(\et_i(b_1\otimes b_2))=\Hb(b_1\otimes
b_2)+\dfrac{\delta_{i0}}{d}
\begin{cases}
1 & \text{if $\vp_0(b_1)\ge\ve_0(b_2)$} \\
-1 & \text{if $\vp_0(b_1)<\ve_0(b_2)$.}
\end{cases}
\end{equation}
The combinatorial $R$-matrices satisfy the Yang-Baxter equation,
which asserts the equality of the two ways to compute the
isomorphism $B_1\otimes B_2\otimes B_3\rightarrow B_3 \otimes
B_2\otimes B_1$ by switching adjacent factors using maps of the form
$R_{B_i,B_j}$.
\item Each $B\in\CC$ has a $U_q(\gf)$-equivariant grading
$\Db_B:B\rightarrow\frac{1}{d}\Z$ called the coenergy function for
$B$, which is well-defined up to a global additive constant by the
following rules \cite{OSS2} \cite{HKOTT}. If $B\in\CC$ is a KR
crystal then the grading $\Db_B$ is prescribed by \cite{HKOTY}
\cite{HKOTT}. If $B_1,B_2\in\CC$ are such that $\Db_{B_1}$ and
$\Db_{B_2}$ have been defined, let $b_1,b_2,b_1',b_2'$ be as in
\eqref{E:H}. Then
\begin{equation} \label{E:D2}
\Db_{B_1\otimes B_2}(b_1 \otimes b_2) = \Hb(b_1\otimes b_2) +
\Db_{B_1}(b_1) + \Db_{B_2}(b_2').
\end{equation}
\end{enumerate}
Suppose $B_1,B_2,\dotsc,B_m\in \CC$. By induction one may prove the
following. Let $b_i\in B_i$ for $1\le i\le m$,
$B=B_1\otimes\dotsm\otimes B_m$, and $b=b_1\otimes\dotsm\otimes
b_m$. Then one may parenthesize the $m$-fold tensor product in any
way, iterate the pairwise construction of \eqref{E:D2}, and the
resulting coenergy function is always given by \cite{HKOTT}
\cite{OSS2}
\begin{equation}\label{E:NYD}
\Db_B(b)=\sum_{1\le i<j\le m}\Hb(b_i\otimes b_{j}^{(i+1)})
+\sum_{j=1}^m \Db_{B_j}(b_j^{(1)})
\end{equation}
where for $1\le i\le j\le m$ the vertices $b_j^{(i)}\in B_j$ are
determined by the affine crystal isomorphisms
\begin{equation} \label{E:RNY}
\begin{split}
B_i \otimes \dotsm \otimes B_{j-1} \otimes B_j &\rightarrow B_j
\otimes B_i \otimes B_{i+1} \otimes\dotsm\otimes B_{j-1}
\\
b_i\otimes \dotsm\otimes b_{j-1}\otimes b_j&\rightarrow b_j^{(i)}
\otimes b_i' \otimes \dotsm \otimes b_{j-1}'%
\end{split}
\end{equation}
given by compositions of combinatorial $R$-matrices acting at
adjacent tensor factors. If each $B_i$ is the same crystal $B$ then
\begin{equation} \label{E:Dhomog}
  \Db_{B^{\otimes m}}(b_1\otimes\dotsm\otimes b_m) =
  \sum_{j=1}^{m-1} (m-j)\Hb(b_j\otimes b_{j+1}) +
  m \Db_B(b_1).
\end{equation}

We now change notation. Up to now the subscript in $B_i$ had merely
indicated the position of a generic tensor factor. {}From now on let
$B_s = B^{1,s}$. Let $B_\mu\in \CC$ be the affine crystal graph
\begin{equation} \label{E:Bmu}
  B_\mu = B_{\mu_1} \otimes B_{\mu_2} \otimes \dotsm \otimes
  B_{\mu_m}
\end{equation}
where $\mu=(\mu_1,\mu_2,\dotsc,\mu_m)\in\Z_{>0}^m$.

For $\la\in P^+$ let $\F_{\la,\mu}$ denote the set of highest weight
vertices in $B_\mu$ of weight $\la$. The 1-d sums for $\CC$ are the
graded tensor product multiplicities defined by
\begin{equation}
\label{E:1dsum}
 \Xb_{\la,\mu}(q) = \sum_{b\in \F_{\la,\mu}} q^{\Db_{B_\mu}(b)}.
\end{equation}

The papers \cite{HKOTY} \cite{HKOTT} introduced another remarkable
expression $M$ called the fermionic formula and conjectured that
$X=M$. The $M$ formula, a sum of products of $q$-binomial
coefficients, arises from the Bethe Ansatz and exhibits
quasiparticle behavior of interest to physicists. We shall not
pursue the $M$ formula here.

\begin{remark} To go from the coenergy functions defined here to the energy
functions of \cite{HKOTY} \cite{HKOTT}, set $d=1$ in the formulae
and take the negative. The factor $d$ is included in our definitions
to make the statement of Theorem \ref{T:X=K} smoother. For the
coenergy to be completely well-defined the local coenergy functions
and the coenergy functions for KR crystals must be normalized. This
is done in section \ref{sub_sec_en} for the cases under
consideration in this paper.
\end{remark}

\subsection{Stable 1-d sums}
Fix a nonexceptional affine algebra $\gaf_n$ together with a
distinguished simple Lie subalgebra $\gf_n$. Let $\Par$ be the set
of partitions and $\Par_m$ those with at most $m$ parts.

\begin{lemma}\label{L:stable} \cite{SZ} \footnote{For technical
reasons this result was stated for the fermionic formula $M$ in
\cite{SZ}.} Fix $m\in\Z_{>0}$, $\la,\mu\in\Par_m$, and a
nonexceptional family of affine algebras $\{\gaf_n\}$. Then for any
$n\ge m$, the 1-d sum $X_{\la,\mu}^{\gaf_n}(q)$ is well-defined and
gives the same polynomial, called a stable 1-d sum. Moreover, all
the families of nonexceptional affine algebras yield only four
distinct kinds of stable 1-d sums, one for each of the partitions
$\vn, (1), (2), (1,1)$ of size at most 2. They are grouped as
follows, based on the attachment of the zero node to the classical
Dynkin diagram of $\gf_n$.
\begin{itemize}
\item Kind $\vn$: $A_{n-1}^{(1)}$: attach $0$ to $1$ and $n-1$ by single
bonds, in the diagram $A_{n-1}$.
\item Kind $(1)$: $D_{n+1}^{(2)}$, $A_{2n}^{(2)\dagger}$, $B_n^{(1)\dagger}$:
attach $0$ to $1$ with a double bond pointing towards $0$, in the
diagrams $B_n$, $C_n$, and $D_n$ respectively.
\item Kind $(2)$: $A_{2n}^{(2)}$, $C_n^{(1)}$,
$A_{2n-1}^{(2)\dagger}$: attach $0$ to $1$ with a double bond
pointing towards $1$, in the diagrams $B_n$, $C_n,$ and $D_n$
respectively.
\item Kind $(1,1)$: $B_n^{(1)}$, $A_{2n-1}^{(2)}$, $D_n^{(1)}$:
Attach $0$ to the vertex $2$ by a single bond, in the diagrams
$B_n$, $C_n,$ and $D_n$ respectively.
\end{itemize}
\end{lemma}
For $\la,\mu\in\Par_m$ and $n\ge m$ we define the stable 1-d sums
\begin{equation} \label{E:stableXdef}
  \Xb_{\la,\mu}^\dia(q) = \Xb_{\la,\mu}^{\gaf_n}(q)
\end{equation}
with $\dia\in\{\vn,(1),(2),(1,1)\}$ and $\gaf_n$ related as above.
The stable 1-d sums of type $A_{n-1}^{(1)}$ satisfy
\begin{equation} \label{E:XKnull}
  \Xb_{\la,\mu}^\vn(q) = \Kb_{\la,\mu}(q) := q^{||\mu||} K_{\la,\mu}(q^{-1})
\end{equation}
where $K_{\la,\mu}(q)$ is the Kostka-Foulkes polynomial
\cite[III.6]{mac}, $\Kb_{\la,\mu}(q)$ is the ``cocharge"
Kostka-Foulkes polynomial, and
\begin{equation} \label{E:mubarbar}
||\mu||=\sum_i (i-1)\mu_i.
\end{equation}
The equality \eqref{E:XKnull} of type $A_{n-1}^{(1)}$ 1-d sums (of
tensor products of symmetric or exterior powers) with the Kostka
Foulkes polynomials, was proved by Nakayashiki and Yamada \cite{NY}
using the combinatorial characterization of Kostka-Foulkes
polynomials due to Lascoux and Sch\"utzenberger \cite{LSc1}.

\subsection{$X=K$}
Let $\Par^\dia$ be the set of partitions that can be tiled with the
partition $\dia$. Then $\Par^\vn$ contains just the empty partition,
$\Par^{(1)}=\Par$ is the set of all partitions, $\Par^{(2)}$
consists of the partitions with even rows and $\Par^{(1,1)}$ the
partitions with even columns. Write $\Par^\dia_m = \Par^\dia \cap
\Par_m$. Let
\begin{equation}\label{E:mubar}
|\mu|=\sum_{i=1}^m \mu_i\qquad\text{for any $\mu\in\Z^m$.}
\end{equation}
The $K$-polynomials \cite{SZ} are defined by
\begin{equation} \label{E:Kdef}
\Kb_{\la,\mu}^\dia(q) = q^{\frac{|\mu|-|\la|}{2}}
  \sum_{\nu\in\Par_m} \Kb_{\nu,\mu}(q) \sum_{\gamma\in \Par_m^\dia}
  c^\nu_{\la\gamma}
\end{equation}
where $c^\nu_{\la\gamma}$ is the Littlewood-Richardson coefficient
\cite[(I.5.2)]{mac}. They are just (up to degree shift) nonnegative
integer combinations of cocharge Kostka-Foulkes polynomials.

Our first main theorem is the $X=K$ theorem, which was conjectured
in \cite{SZ}.

\begin{theorem} \label{T:X=K} ($X=K$) Let $\la,\mu\in\Par_m$. Then for every
$\dia\in \{\vn,(1),(2),(1,1)\}$,
\begin{equation} \label{E:X=K}
  \Xb_{\la,\mu}^\dia(q) = \Kb_{\la,\mu}^\dia(q).
\end{equation}

\end{theorem}

The case $\dia=\vn$ is trivial. The cases $\dia=(1)$ and $\dia=(2)$
were proved in \cite{sh} using the virtual crystal construction of
\cite{OSS2}, which embeds affine crystals of types $D_{n+1}^{(2)}$,
$A_{2n}^{(2)}$, and $C_n^{(1)}$ into those of type $A_{2n-1}^{(1)}$.
We establish the remaining case $\dia=(1,1)$ in section
\ref{sec_proof1}. Since such an affine crystal embedding into type
$A^{(1)}$ does not exist for types $B_n^{(1)}$, $A_{2n-1}^{(2)}$, or
$D_n^{(1)}$, we must use an entirely different approach to relate
the 1-d sums $\Xb_{\la,\mu}^{(1,1)}(q)$ to those of type $A^{(1)}$.

\subsection{Lusztig $q$-analogues}
Let $\rho=\frac{1}{2}\sum_{\alpha\in R^+} \alpha=\sum_{i\in\I}
\omega_i$. Fix a function $\lwt:R^+\rightarrow\Z$ that is constant
on the orbits of the Weyl group $W$. For $\la,\mu\in P^+$, define
the polynomial
\begin{equation} \label{E:Lusztig}
  \KL_{\la,\mu}^{\gf,L}(q) = \sum_{w\in W} (-1)^w [e^{w(\la+\rho)-(\mu+\rho)}] \prod_{\alpha\in R^+}
  \dfrac{1}{1-q^{\lwt(\alpha)} e^\alpha}
\end{equation}
where $[e^\beta]f$ denotes the coefficient of $e^\beta$ in $f\in
\Z[P]$. It follows from the Weyl character formula that
$\KL_{\la,\mu}^{\gf,L}(1)$ is the multiplicity of the weight $\mu$
in $V(\la)$. When $\lwt(\alpha)=1$ for all $\alpha\in R^+$,
$\KL_{\la,\mu}^{\gf,L}(q)=\KL_{\la,\mu}^\gf(q)$ is Lusztig's
$q$-analogue of weight multiplicity \cite{Lu}.

Lusztig's $q$-analogues are certain affine Kazhdan-Lusztig
polynomials. Let $\tW\cong P \rtimes W$ be the extended affine Weyl
group, realized as a subgroup of isometries of the weight lattice
$P$. For $\la\in P$ let $t_\la\in\tW$ be translation by $\la$. For
$\la\in P^+$ let $w_\la$ be the element of maximal length in the
double coset $W t_\la W\subset \tW$. Then
\begin{equation} \label{E:qmultKL}
  \KL_{\la,\mu}^{\gf}(q) = q^{(\ell(w_\la)-\ell(w_\mu))/2} P_{w_\mu,w_\la}^{\gaf}(q^{-1})
\end{equation}
where $P_{x,y}^{\gaf}(q)$ is the Kazhdan-Lusztig polynomial
\cite{Lu}.

\subsection{Stable KL polynomials} We consider a subfamily
of Lusztig $q$-analogues that first studied in \cite{lec}
\cite{lec4}. Let $\la,\mu\in\Par_m$ with $n\ge m$, regarded as
dominant weights for $\gf\in\{A_{n-1},B_n,C_n,D_n\}$. Define the
stable KL polynomial
\begin{equation} \label{E:stableKL}
\KLS_{\la,\mu}^{\gf_n,L}(q) = \sum_{w\in S_n} (-1)^{\ell(w)} [
e^{w(\la+\rho)-(\mu+\rho)} ] \prod_{\alpha\in R^+(\gf_n)} (1-
q^{L(\alpha)} e^\alpha)^{-1}.
\end{equation}
This is an expression for $\gf_n$, but the sum runs over the
parabolic subgroup of $W$ generated by $s_i$ for $1\le i\le n-1$;
this is a copy of the symmetric group $S_n=W(A_{n-1})$, the Weyl
group of type $A_{n-1}$. The following stability phenomenon, which
justifies the above nomenclature, was observed in \cite[Theorem
5.1.5]{lec}.

\begin{proposition} \label{P:stableKL} With $\la,\mu$ as above,
we have
\begin{equation} \label{E:KL=stableKL}
\KL_{\la+(k^n),\mu+(k^n)}^{\gf_n,L}(q) = \KLS_{\la,\mu}^{\gf_n,L}(q)
\end{equation}
provided that $k \ge (|\la|-|\mu|)/2$. In particular the left hand
side of \eqref{E:KL=stableKL} does not change for $k$ sufficiently
large. Moreover for $\gf_n=A_{n-1}$ every Lusztig $q$-analogue is a
stable KL polynomial and Kostka-Foulkes polynomial
\begin{equation} \label{E:stableKLA=KF}%
\KLS_{\la,\mu}^{A_{n-1}}(q) = \KL_{\la,\mu}^{A_{n-1}}(q) =
K_{\la,\mu}(q).
\end{equation}
\end{proposition}
\begin{proof} It is well known that $\la\ge w\la$ for all
$w\in W$ and $\la\in P^+$, where
\begin{equation} \label{E:order}
\la\ge\mu\qquad\text{if and only if}\qquad \la-\mu\in Q^+ =
\bigoplus_{i\in I} \Z_{\ge0}\alpha_i.
\end{equation}
In particular $w(\la+\rho)-(\mu+\rho)\le \la-\mu$. From this it
follows that
\begin{equation} \label{E:KLzero}
\KL_{\la,\mu}^{\gf_n,L}(q)=0 \qquad\text{unless $\la\ge\mu$}
\end{equation}
for if $\la\not\ge\mu$ then every summand is zero. Note that the
definition \eqref{E:stableKL} makes sense for pairs of elements in
the set
\begin{equation} \label{E:DWdef}
\DW = \{ \la\in\Z^n\mid \la_1\ge\la_2\ge\dotsm\ge\la_n \}
\end{equation}
where the parts of $\la$ may be negative. Since every $w\in S_n$
fixes all multiples of $(1^n)$,
\begin{equation} \label{E:KLshift}
\KLS_{\la+(k^n),\mu+(k^n)}^{\gf_n,L}(q)=\KLS_{\la,\mu}^{\gf_n,L}(q)\qquad\text{for
all $k\in\Z$.}
\end{equation}
Let us assume that $k \ge (|\la|-|\mu|)/2$, recalling the definition
\eqref{E:mubar}.

In the case $\gf_n=A_{n-1}$, since $\la-\mu\in Q^+$ and
$|\alpha_i|=0$ for all $1\le i\le n-1$, we have $|\la|=|\mu|$. Thus
all Lusztig $q$-analogues are stable KL polynomials and
\eqref{E:KL=stableKL} holds since $W(A_{n-1})=S_n$. The equality
\eqref{E:stableKLA=KF} may be taken as the definition of the
Kostka-Foulkes polynomials.

Now let $\gf_n\in\{B_n,C_n,D_n\}$. Suppose that $w\in W\backslash
S_n$. It suffices to show that the $w$-th summand in
$\KL_{\la+(k^n),\mu+(k^n)}^{\gf_n,L}(q)$ is zero. Suppose not. Then
$\beta=w(\la+(k^n)+\rho)-(\mu+(k^n)+\rho)\in Q^+$. In particular
$|\beta|=\inner{\beta}{(1^n)}\ge0$ since $(1^n)$ is a positive
integer multiple of a fundamental weight. Since $w\not\in S_n$,
$w(k^n)$ must have some coordinate $-k$ so that $|w(k^n)-(k^n)| \le
-2k$. We also have $|w(\rho)-\rho|<0$ since $\rho$ has nonnegative
coordinates and $w(\rho)$ has a strictly negative coordinate. We
obtain the contradiction
\begin{equation*}
0 \le |w(\la+(k^n)+\rho)-(\mu+(k^n)+\rho)| \le
|w(k^n)-(k^n)|+|w(\rho)-\rho|+|w(\la)-\mu| < -2k + |\la|-|\mu| \le
0.
\end{equation*}
\end{proof}

\subsection{$K=\KLS$}
We now state the $K=\KLS$ theorem, which asserts the equality of the
$K$-polynomials and the stable KL polynomials. Let
$\la,\mu\in\Par_m$ with $n\ge m$. Define the partitions
\begin{equation}\label{E:hat}
\begin{split}
\widehat{\la}&=(M-\la_m,\dotsc,M-\la_1) \\
\widehat{\mu}&=(M-\mu_m,\dotsc,M-\mu_1)
\end{split}
\qquad \qquad\text{where $M\ge \max(\la_1,\mu_1)+(|\mu|-|\la|)/2$.}
\end{equation}

\begin{theorem} \label{T:K=KL} ($K=\KLS$) We have
\begin{equation} \label{E:K=KL}
\begin{split}
  \Kb_{\la,\mu}^\dia(q) &= q^{||\mu||+|\mu|-|\la|}\,\,
  \KLS_{\lah,\muh}^{\gf_n,L}(q^{-1}) \\
  &= q^{||\mu||+|\mu|-|\la|}\,\,
  \KL_{\lah,\muh}^{\gf_n,L}(q^{-1})
\end{split}
\end{equation}
where $\dia=\{\vn,(1),(2),(1,1)\}$ respectively correspond to
$\gf_n=\{A_{n-1},B_n,C_n,D_n\}$ and $L(\alpha)=1$ for all $\alpha\in
R^+$ except for $\gf_n=B_n$, in which case $L(\alpha)=\frac{1}{2}$
for short roots.
\end{theorem}

The choice of $M$ guarantees that the KL polynomials
$\KL_{\lah,\muh}^{\gf_n,L}(q)$ are stable; the second equality in
\eqref{E:K=KL} holds by Proposition \ref{P:stableKL}. Theorem
\ref{T:K=KL} proved in section \ref{S:KKL} by a type-independent
argument.

For $\dia=\vn$, $K=\KLS$ reduces to the definition of
$\Kb_{\la,\mu}(q)$ in \eqref{E:XKnull} due to a duality satisfied by
the Kostka-Foulkes polynomials (see \cite{LS} or \eqref{E:KLAdual}):
\begin{equation} \label{E:KFdual}
  K_{\la,\mu}(q) = K_{\lah,\muh}(q).
\end{equation}

For $\dia=(2)$ and $\dia=(1,1)$, $K=KL$ was already proved in
\cite{lec} \cite{lec4}.

For $\dia=(1)$, in \cite{lec4} it was observed that \eqref{E:K=KL}
is false if one takes $L(\alpha)=1$ for all $\alpha\in R^+$, that
is, when the right hand side is the usual Lusztig $q$-analogue for
type $B_n$.

Combining $X=K$ and $K=\KLS$ (Theorems \ref{T:X=K} and \ref{T:K=KL})
we obtain

\begin{corollary} \label{C:X=KL} ($X=\KLS$ Theorem) With notation as
in Theorem \ref{T:K=KL},
\begin{equation} \label{E:X=KL}
  \Xb_{\la,\mu}^\dia(q) = q^{||\mu||+|\mu|-|\la|}\,\,
  \KL_{\lah,\muh}^{\gf_n,L}(q^{-1})
\end{equation}
\end{corollary}

\section{Explicit description of crystals and energy functions}
In this section we give details for the 1-d sums involved in our
proof of $X=K$ for $\dia=(1,1)$.

\subsection{Notation for KR crystals}
\label{S:KR}

Let $\mu \in\Z_{>0}^m$. By abuse of notation, denote by
$B_\mu^{\gf_n}$ the affine crystal graph given by the tensor product
defined in \eqref{E:Bmu} for the affine root systems $\gaf_n\in
\{A_{n-1}^{(1)},A_{2n-1}^{(2)},D_n^{(1)\dagger}\}$, which have respective
simple Lie subalgebras $\gf_n\in\{A_{n-1},C_n,D_n^\dagger\}$; here $D_n^{(1)\dagger}$
is the affine root system $D_n^{(1)}$ except that $n$ is regarded as the affine node,
so that $D_n^\dagger$ is its classical sub root system. We have
\cite{KMN2} \cite{KKM}
\begin{equation} \label{E:KRconnected}
B_s^\gf\cong B^\gf(s\omega_1)\qquad\text{as $\gf$-crystals for
$\gaf_n\in\{A_{n-1}^{(1)},A_{2n-1}^{(2)},D_n^{(1)\dagger}\}$.}
\end{equation}
{}From now on we shall use $D_n^\dagger$ and $D_n^{(1)\dagger}$ instead
of $D_n$ and $D_n^{(1)}$ unless specifically indicated otherwise.

\subsection{$A_{n-1}$, $C_n$ and $D_n$ crystals}
The crystals of type $A_{n-1},C_n,D_n^\dagger$ are colored directed
graphs with colors in the respective sets $\J^A=\{1,2,\dotsc,n-1\}$,
$\J=\{1,2,\dotsc,n\}$, and $\J^\dagger=\{0,1,2,\dotsc,n-1\}$. The
$U_q(\gf)$-modules associated to the partition
$\la=(1,0,\dotsc,0)\in\Par_n$ are the vector representations
$V^{A_{n-1}}(\omega_1)$, $V^{C_n}(\omega_1)$, and
$V^{D_n^\dagger}(\omega_{n-1}^\dagger)$ respectively. Their
respective dimensions are $n,2n,2n$. We call these crystal graphs
$B_1^{A_{n-1}}$, $B_1^{C_n}$, and $B_1^{D_n}$ respectively. This
agrees with the notation of the previous section, except that the
crystal graphs pictured above are $\gf$-crystals as opposed to
$\gaf$-crystals. The crystal graphs $B_1^\gf$ are depicted in Figure
\ref{F:B1}. The notation $B_1^{D_n}$ is consistent with the usual
Dynkin labeling $D_n$, but the pictured crystal graph uses the
conventions of $D_n^\dagger$. We use the set of symbols
$\{1,2,\dotsc,n,\bn,\dotsc,\bbb,\ba\}$, with weights given by
$\wt(i)=\ve_i$ and $\wt(\overline{i})=-\ve_i$ for $1\le i\le n$.
\begin{figure}
\begin{align*}
&B_1^{A_{n-1}} &\quad& \fb{1}\overset{1}{\longrightarrow}
\fb{2}\overset{2}\longrightarrow\dotsm\longrightarrow
\fbs{n\!\!-\!\!1} \overset {n-1}{\longrightarrow}\fb{n} \\[2mm]
&B_1^{C_n}&\quad&
\fb{1}\overset{1}{\longrightarrow}\fb{2}\overset{2}\longrightarrow\dotsm\longrightarrow
\fbs{n\!\!-\!\!1} \overset
{n-1}{\longrightarrow}\fb{n}\overset{n}{\longrightarrow}\fb{\bn}\overset{n-1}%
{\longrightarrow}\fbs{\overline{n\!\!-\!\!1}}\overset{n-2}{\longrightarrow}\dotsm
\longrightarrow\fb{\bbb} \overset{1}{\longrightarrow}\fb{\ba}
\\[2mm]
&B_1^{D_n}&\quad&
\fb{\bn}\overset{n-1}{\longrightarrow}\fbs{\bnm}\overset
{n-2}{\longrightarrow}\dotsm\overset{3}{\longrightarrow}\fb{\overline{3}}%
\overset{2}{\longrightarrow}%
\begin{tabular}
[c]{c}%
\ $\fb{\overline{1}}$ \ \ \\[-1mm]
\ \ $\overset{1\ \ }{\nearrow}$ $\ \ \ \overset{\ \ 0}{\searrow}$ \ \ \ \\
$\fb{\overline{2}}\ \ \ \ \ \ \ \ \ \ \ \ \ \fb{2}$\\
\ $\underset{0\ \ }{\searrow}$ \ \ \ $\underset{\ \ 1}{\nearrow}$ \ \ \
\\[-1mm]
$\fb{1}$ \
\end{tabular}
\overset{2}{\longrightarrow}\fb{3}\overset{3}{\longrightarrow}\dotsm\overset
{n-2}{\longrightarrow}\fbs{n\!\!-\!\!1}\overset{n-1}{\longrightarrow}\fb{n} \label{vect_D}%
\end{align*}
\caption{\label{F:B1} Crystal graphs $B_1^\gf$}
\end{figure}

We define a partial order $\le^\gf$ on each set $B_1^\gf$ based on
reachability.
\begin{equation*}
\begin{split}
(B_1^{A_{n-1}},\le^{A_{n-1}}) &= \{1<2<\dotsm<n\} \\
(B_1^{C_n},\le^{C_n}) &=\{1<2<\dotsm<n<\bn<\dotsm <\overline{2}<\overline{1}\} \\
(B_1^{D_n},\le^{D_n}) &= \{ \bn < \dotsm < \bbb < \begin{matrix} 1 \\[2mm] \ba
\end{matrix}
<  2\dotsm < n\}.
\end{split}
\end{equation*}
The letters $1$ and $\overline{1}$ are not comparable for
$\le^{D_n}$.

For $\gf_n\in\{A_{n-1},C_n,D_n\}$ there is an obvious bijection
between the vertices of $B_{1^s}^\gf = (B_1^\gf)^{\otimes s}$ and
the words of length $s$ in the set of symbols $B_1^\gf$ given by
$b=x_1\otimes x_2\otimes\dotsm\otimes x_l\mapsto w(b):=x_1x_2\dotsm
x_l$. Denote its inverse by $w\mapsto b(w)$.

\subsection{Tableaux of type $A_{n-1},C_n,D_n$}
Let $\gf\in\{A_{n-1},C_n,D_n\}$. Kashiwara and Nakashima \cite{KN}
introduced $\gf$-tableaux, which for $\gf=A_{n-1}$ are the
well-known semistandard tableaux. For $\la\in\Par_n$, the
$\gf$-tableaux of shape $\la$ give a natural labeling of the crystal
$B^\gf(\la)$. They are fillings of the Young diagram associated to
$\la$ which are semistandard for the orders defined on the set of
symbols $B_1^\gf$ and satisfy additional conditions on their rows
and columns detailed in \cite{KN} and \cite{lec2}. Write $\TT^\gf$
for the set of $\gf$-tableaux.

\begin{example} $T_1\in \TT^{C_4}$ and $T_2\in\TT^{D_4}$ where
\[%
\Yboxdim{14pt}%
 T_1 = \young(13\bc,\bc\bc,\bbb\bbb)\qquad T_2 =
\young(\bd\bc1,\bbb\bbb,23).
\]
\end{example}

In the sequel we identify the word $x_1\dotsm x_s$ with the row
tableau
\[
L=%
\begin{tabular}
[c]{|c|c|c|c|c|c|}\hline $\!\!x_s\!\!$ & $\cdot$ & $\cdot$ & $\cdot$
& $\!\!x_2\!\!$ & $\!\!x_1\!\!$ \\ \hline
\end{tabular}
\]%

\begin{proposition}
\label{prop_KN} \cite{KN} For $\gf\in\{A_{n-1},C_n,D_n\}$, the
vertices of $B_s^\gf$ can be labeled by the decreasing words of
length $s$ on the set $B_1^\gf$, that is, by the words $x_1x_2\dotsm
x_s$ with $x_1\ge^\gf x_2 \ge^\gf \dotsm \ge^\gf x_s$.
\end{proposition}

In particular a decreasing word of type $D_n$ cannot contain both
letters $1$ and $\ba$.

\begin{remark} \label{R:AinCD} Let $\gf=C_n,D_n$ and $\mu\in\Z_{>0}^m$ with $n\ge m$.
The subset of elements in $B_\mu^{\gf}$ that involve no barred
letters, can be identified with $B_\mu^{A_{n-1}}$ as
$A_{n-1}$-crystals. So we can write $B_\mu^{A_{n-1}} \subset
B_\mu^\gf$.
\end{remark}

\subsection{Plactic monoids and insertion algorithms for types $A_{n-1},C_n,D_n$}
Let $\gf\in\{A_{n-1},C_n,D_n\}$. Consider the tensor crystal $G^\gf
= \bigoplus_{l\ge0} (B_1^\gf)^{\otimes l}$. Denote by $\sim^\gf$ the
equivalence relation on the vertices of $G^\gf$ defined by $b
\sim^\gf b'$ if $b$ and $b'$ are in the same connected component.
The quotient set $G^\gf/\!\!\sim^\gf$ is a monoid \cite{Lec}
\cite{lec2} \cite{lit} which is a quotient of the free monoid on the
set $B_1^\gf$ by two kinds of plactic relations. The first consists
of relations of length $3$. They reflect the crystal isomorphism
$B(121)\simeq B(112)$. For $\gf=A_{n-1}$ they are the Knuth
relations which define Lascoux-Sch\"utzenberger's plactic monoid
\cite{LS}:
\begin{align*}
  yzx&\equiv yxz &\qquad &\text{for $x\le y <z$ and } \\
  xzy&\equiv zxy &\qquad &\text{for $x <y\le z$.}
\end{align*}
For $\gf=C_n$ we obtain
\[
\left\{
\begin{tabular}
[c]{l}%
$R_1^C:$ $yzx\equiv yxz$ for $x\le^C y<^Cz$ with
$z\neq\overline{x},$ and
$xzy\equiv zxy$ for $x<^Cy\le^C z$ with $z\neq\overline{x}\vspace*{0.2cm}$\\
$R_2^C:y(\overline{x-1})(x-1)\equiv yx\overline{x}$ and $x\overline
{x}y\equiv(\overline{x-1})(x-1)y$ for $1<^Cx\le^C n$ and $x\le^C
y\le^C\overline
{x}$%
\end{tabular}
\right.
\]
and for $\gf=D_n$ we have:%
\[
\left\{
\begin{tabular}
[c]{l}%
$R_1^D:$ If $x\neq\overline{z},$ $yzx\equiv yxz$ for $x\le^D y<^D z$
\ and \ $xzy\equiv zxy$ for $x <^D y\le^D z\vspace*{0.2cm}$\\
$R_2^D:$ If $1<^D x<^D n$ and $\overline{x}\le^D y\le^D x,$
$y(x+1)\overline{(x+1)}\equiv y\overline{x}x$ \ and \ $\overline{x}%
xy\equiv(x+1)(\overline{x+1})y\vspace*{0.2cm}$\\
$R_{3}^D:$ If $x\notin\{1,\overline{1}\},$ $\left\{
\begin{tabular}
[c]{l}%
$1\,x\overline{1}\equiv\overline{x}\,1\overline{1}$\\
$\overline{1}\,x\,1\equiv x\,\overline{1}1$%
\end{tabular}
\right.  \text{ and }\left\{
\begin{tabular}
[c]{l}%
$1\overline{1}\overline{x}\equiv1\overline{x}\overline{1}$\\
$\overline{1}1\overline{x}\equiv\overline{1}\overline{x}1$%
\end{tabular}
\right.  \vspace*{0.2cm}$\\
$R_{4}^D:\left\{
\begin{tabular}
[c]{l}%
$\overline{1}\,11\equiv2\overline{2}\,1$\\
$1\overline{1}\,\overline{1}\equiv2\overline{2}\,\overline{1}$%
\end{tabular}
\right.  \text{ and }\left\{
\begin{tabular}
[c]{l}%
$12\overline{2}\equiv11\overline{1}$\\
$\overline{1}2\overline{2}\equiv\overline{1}\,\,\overline{1}1$%
\end{tabular}
\right.  $%
\end{tabular}
\right.  .
\]
The second kind of relation is called a ``contraction relation''.
This kind of relation does not preserve the length of words. These
relations reflect the crystals isomorphisms $B(1\dotsm p\overline
{p})\simeq B(1\dotsm p-1)$ with $p\in\{1,\dotsc,n\}$ for $\gf=C_n$
and $B(\bn\dotsm\overline{q}q)\simeq B(\bn%
\dotsm\overline{q+1})$ with $q\in\{1,\dotsc,n\}$ for $\gf=D_n$. The
reader is referred to \cite{Lec} and \cite{lec2} for a complete
description of these relations.

For each $\gf$ there is an insertion scheme that is compatible with
the above plactic relations for $\gf$-tableaux \cite{Lec}
\cite{lec2}. We shall call this $\gf$-insertion. $A_{n-1}$-insertion
is the well-known Robinson-Schensted insertion algorithm on
semistandard tableaux \cite{Ful}. Consider a word $w=x_1\dotsm x_s$
with $x_i\in B_1^\gf$ such that $\mathrm{b}(w)$ belongs to a
connected component of $G^\gf$ isomorphic to $B^\gf(\la)$.
$\gf$-insertion computes the unique $\gf$-tableau $P^\gf(w)$
occurring in $B^\gf(\la)$ such that $\mathrm{b}(w)\sim^\gf
P^\gf(w)$.

To define $\gf$-insertion it suffices to describe the insertion
denoted $x \overset{\gf}{\rightarrow} T$ of the letter $x$ into the
tableau $T$. This insertion is obtained directly from the plactic
relations. For any word $w=x_1\dotsm x_s$, $P^\gf(w)$ is defined
recursively by
\[
P^\gf(x_1\dotsm x_i)=(x_i \overset{\gf}{\rightarrow} P^\gf(x_1\dotsm
x_{i-1}))
\]
where the $\gf$-insertion of the empty word results in the empty
tableau. For complete details on $\gf$-insertion see \cite{Lec}
\cite{lec2}.

For $\delta\in\Z_{>0}^m$ let $F_\delta^{A_{n-1}}$ be the set of
$A_{n-1}$-highest weight vectors in $B_\delta^{A_{n-1}}$.

\begin{lemma}
\label{lemm_B(nu)} For any $\delta\in\Z_{>0}^m$ with $m<n$ and
$\nu\in\Par_m$, there is a bijection $\Psi$ from the set of
$A_{n-1}$-components of $B_\delta^{A_{n-1}}$, to the set of
$D_n$-components of $B_\delta^{D_n}$ whose highest weight is a
partition $\nu$ such that $|\nu|=|\delta|$. It is defined by
\begin{equation} \label{E:ADcomp}
  \Psi(B^{A_{n-1}}(b)) = B^{D_n}(b)
\end{equation}
for any $b\in F_\delta^{A_{n-1}}$ where $B^{A_{n-1}}(b)$ is the
$A_{n-1}$-component of $b$ in $B_\delta^{A_{n-1}}$ and $B^{D_n}(b)$
is the $D_n$-component of $b$ in $B_\delta^{D_n}$.
\end{lemma}

\begin{remark} \label{R:T}
\begin{enumerate}
\item $\TT^{D_{n-1}}\subset\TT^{D_n}$.
\item For $T\in\TT^{D_{n-1}}$ and
$b=x_1\dotsm x_s\in B_s^{D_n}$, no contraction occurs during the
insertions
\[%
x_l\overset{D_n}{\mathbf{\longrightarrow}}(x_2\overset
{D_n}{\longrightarrow}(\dotsm(x_1\overset{D_n}{\longrightarrow}T))).%
\]
\item Let $b\in B_\delta^{D_n}$. Then $B^{D_n}(b)$ is in the
image of the map $\Psi$ of \eqref{E:ADcomp} if and only if there is
no contraction during the insertion procedure $P^{D_n}(b_1b_2\dotsm
b_m)$.
\end{enumerate}
\end{remark}

\subsection{Left splitting embeddings}

\begin{lemma}
\label{ele_cry_iso} For any nonnegative integers $k,l$ and for
$\gf=A_{n-1},C_n,D_n$ there is a unique $\gf$-crystal isomorphism
\[
\sigma:B_l\otimes B_k\rightarrow B_k\otimes B_l
\]
\end{lemma}
\begin{proof} For $\gf=A_{n-1},C_n,D_n$ this follows respectively from the
uniqueness of the combinatorial $R$-matrix for $\gaf=A_{n-1}^{(1)},
A_{2n-1}^{(2)}, D_n^{(1)}$, acting on the affine crystal $B_l
\otimes B_k$.
\end{proof}

\begin{remark} \label{R:tensor2}
Lemma \ref{ele_cry_iso} follows more simply from the fact that in
the above cases, $B_l \otimes B_k$ is multiplicity-free. This is
seen from the following list of highest weight vertices of
$B_l\otimes B_k$ computed by Lemma \ref{lem_plu_hp}.
\begin{align*}
\label{high_n=2} v_{l,k;b}^{A_{n-1}}&=1^l\otimes
2^b1^{k-b}&\qquad&\text{$0\le b\le \min(k,l)$} \\
v_{l,k;a,b}^{C_n}&=1^l\otimes \overline{1}^a 2^b1^{k-a-b}%
&\qquad&\text{$a,b\ge0$, $a+b\le \min(k,l)$} \\
v_{l,k;a,b}^{D_n}&=\bn^l\otimes n^a \overline{n\!-\!1}^b
\bn^{k-a-b}&\qquad&\text{$a,b\ge0$, $a+b\le\min(k,l)$.}
\end{align*}
\end{remark}

\begin{remark} The crystal isomorphism $\sigma^\gf$ can be computed
using insertion for $\gf$-tableaux \cite{lec2} \cite{Lec} \cite{NY}
\cite{Ok} \cite{Ok2}.
\end{remark}

Let $\gf=A_{n-1},C_n,D_n$. Let $b=b_1\otimes\dotsm\otimes b_m\in
B_\delta$ and $b_1=x_1x_2\dotsm x_{\delta_1}$. The left splitting
operation is defined by
\begin{equation} \label{E:S}
\begin{split}
S:B_\delta &\rightarrow B_{\delta'} \\
b&\mapsto x_1\otimes b_1'\otimes b_2\otimes\dotsm\otimes b_m
\end{split}
\end{equation}
where $\delta'=(1,\delta_1-1,\delta_2,\dotsc,\delta_m)$ and
$b_1'=x_2\dotsm x_{\delta_1}$.

Suppose $\delta_i=1$ for some particular $i$. Write $b_i=x_i$ where
$x_i$ is a letter. Let
\begin{equation}
\begin{split}
\sigma_i:B_\delta&\hookrightarrow B_{\sigma_i(\delta)} \\
b_1\otimes\dotsm\otimes x_i\otimes b_{i+1}\otimes\cdot
\cdot\cdot\otimes b_m&\mapsto b_1\otimes\dotsm\otimes
b_i^{\prime}\otimes x_i^{\prime}\otimes\dotsm\otimes b_m
 \label{sigmai}%
\end{split}
\end{equation}
where $\sigma_i(\delta)$ is obtained by switching the $i$-th and the
$i+1$-th parts of $\delta$ and $\sigma(x_i\otimes
b_{i+1})=b_i'\otimes x_i'$ is the crystal isomorphism of Lemma
\ref{ele_cry_iso}.

By composing left splitting and switching operations one may define
$\gf$-crystal embeddings \cite{OSS}
\begin{equation} \label{E:Sdef}%
S_{\delta}:B_\delta\hookrightarrow B_{1^{|\delta|}}
\end{equation}
We define $S_\delta$ by descending induction on the number of parts
of $\delta$ equal to $1$, and then by the minimum index $k$ such
that $\delta_k > 1$. $S_\delta$ is defined to be the identity if
every part of $\delta$ is $1$. Otherwise some part of $\delta$ is
greater than $1$. Let $k$ be as above. If $k=1$ then we first split
the left tensor factor and induct: $S_\delta = S_{\delta'} \circ S$
with the above notation for $S:B_\delta \rightarrow B_{\delta'}$. If
$k>1$, let $s=\delta_k$. We move this part closer to the front using
the isomorphism $\psi_s:B_1 \otimes B_s \rightarrow B_s \otimes B_1$
acting at the $(k-1)$-th and $k$-th tensor positions, and then
induct. That is, writing $s_{k-1} \delta$ for $\delta$ with its
$(k-1)$-th part $1=\delta_{k-1}$ and $k$-th part $s=\delta_k$
interchanged, define $S_\delta = S_{s_{k-1}\delta} \circ \psi_s$
where it is understood that $\psi_s$ acts at the $(k-1)$-th and
$k$-th tensor factors.

\begin{remark} One may take a shortcut:
instead of applying a splitting map $S:B_s \rightarrow B_1 \otimes
B_{s-1}$ on the leftmost tensor factor, one may split
$S:B_s\rightarrow B_1^{\otimes s}$ all at once, chopping a word of
length $s$ in $B_s$ into its constituent letters.
\end{remark}

\begin{example} \label{ex_spli} Let $x^{\otimes p}=x\otimes
x\otimes\dotsm\otimes x$ denote the $p$-th tensor power of $x$.
\begin{enumerate}
\item Computing the map $S_{(l,k)}^C$ on the element $v_{l,k;a,b}^C$
of Remark \ref{R:tensor2}, the following elements occur as
intermediate values.
\begin{gather*}
1^l \otimes \ba^a 2^b 1^{k-a-b} \qquad 1^{\otimes l} \otimes \ba^a
2^b 1^{k-a-b} \qquad 1^{\otimes l-a} \otimes 2^b 1^{k-b}\otimes
\ba^{\otimes a} \qquad 1^{\otimes l-a-b}\otimes 1^k\otimes
2^{\otimes b}\otimes
\ba^{\otimes a}\\
1^k\otimes1^{\otimes l-a-b}\otimes2^{\otimes b}\otimes \ba^{\otimes
a}\qquad S_{(l,k)}^C(v_{l,k;a,b})=1^{\otimes
k+l-a-b}\otimes2^{\otimes b}\otimes\ba^{\otimes a}%
\end{gather*}
\item For $b=1^4\otimes\ba21\otimes\ba\overline{2}1\in B_{(4,3,3)}^D$
we obtain
$S_{(4,3,3)}^D(b)=1^{\otimes5}\otimes2\otimes\ba^{\otimes2}
\otimes\bbb\otimes1$ in $B_{1^{10}}^D$.
\end{enumerate}
\end{example}

\begin{remark} \label{R:SAinCD}
Let $\gf=C_n,D_n$. With $B_\delta^{A_{n-1}}\subset B_\delta^\gf$ as
in Remark \ref{R:AinCD}, the restriction of the map $S_\delta^\gf$
to $B_\delta^{A_{n-1}}$ agrees with the $A_{n-1}$-crystal embedding
$S_\delta^{A_{n-1}}$.
\end{remark}

\begin{remark} \label{R:indep} One need not split and
permute in exactly the order prescribed in the definition of
$S_\delta$. One may compute $S_\delta$ by applying any sequence of
operations going from $S_\delta$ to $S_{1^{|\delta|}}$ comprised of
splitting the leftmost tensor factor and permuting factors by
combinatorial $R$-matrices.
\end{remark}

\subsection{Coenergy functions}
\label{sub_sec_en} In this section we follow \cite{KMN} \cite{NY}
\cite{HKOTT} \cite{OSS2}.

>From now on, we write $\Db$ and $\Hb$ (resp. $\Db^A$ and $\Hb^A$)
for the coenergy functions associated to affine crystals of type
$A_{2n-1}^{(2)}$ (resp. $A_{n-1}^{(1)}$). By definition, to specify
the coenergy function $\Db$ for the affine crystal $B_\mu^{C_n}$ of
type $A_{2n-1}^{(2)}$, it suffices to specify $\Db_{B_k}$ and
$\Hb=\Hb_{B_l,B_k}$ for the affine crystals $B_l$ and $B_k$ of type
$A_{2n-1}^{(2)}$ for $l\ge k$, and similarly for $\Db^A$.

By \eqref{E:KRconnected} for $\gaf_n=A_{2n-1}^{(2)}$ and
$\gaf_n=A_{n-1}^{(1)}$, $\Db_{B_k}$ and $\Db_{B_k}^A$ are constant
functions. We normalize them to have the value $0$:
\begin{equation} \label{E:DKRnorm}
\Db_{B_k} \equiv 0\qquad \Db_{B_k}^A \equiv 0.
\end{equation}

To normalize $\Hb$ and $\Hb^A$ on $B_l\otimes B_k$ it suffices to
specify a single function value. We set
\begin{equation} \label{E:Hnorm}
  \Hb(1^l \otimes 1^k) = 0 \qquad
  \Hb^A(1^l \otimes 1^k) = 0.
\end{equation}

\begin{remark} \label{R:Hb} Let $\Hb$ denote either $\Hb_{B_l,B_k}$
for $A_{2n-1}^{(2)}$ crystals, or $\Hb_{B_l,B_k}^A$ for
$A_{n-1}^{(1)}$-crystals. Then $\Hb$ can be computed using
respectively $C_n$- or $A_{n-1}$-insertion \cite{Ok} \cite{NY}.
Given a vertex $b=b_1\otimes b_2\in B_l\otimes B_k$, we have
$\Hb(b)=k+l-m$ where $m$ is the size of the first row of the tableau
$P^\gf(b)$ obtained by inserting the row $b_2$ into the row tableau
$b_1$.
\end{remark}

\begin{example}
\label{exa1} For the highest weight vertices $v_{l,k;a,b}^{C_n}\in
B_l^{C_n}\otimes B_k^{C_n}$ in Remark \ref{R:tensor2} we have
\begin{equation*}
\Yboxdim{16pt}%
  P^{C_n}(v_{l,k;a,b}) = \young(1\dotsm\dotsm1,2\dotsm2)
\end{equation*}
which has shape $(l+k-2a-b,b)$. Thus $\Hb(v_{l,k;a,b})=2a+b$.
\end{example}

\begin{example} \label{E:H1}
For $x,y\in B_1^{C_n}$, $\Hb_{B_1,B_1}$ is given by
\begin{equation} \label{H_simple}
\Hb(x\otimes y)=
\begin{cases}
0 & \text{if $x\ge y$} \\
1 & \text{if $x<y$ and $(x,y)\neq(1,\ba)$} \\
2 & \text{if $(x,y)=(1,\ba)$.}
\end{cases}
\end{equation}
For $x,y\in B_1^{A_{n-1}}$ $\Hb^A_{B_1,B_1}$ is given by
\begin{equation} \label{H:A}
\Hb^A(x\otimes y) =
\begin{cases}
0 &\text{if $x\ge y$} \\
1 &\text{otherwise.}
\end{cases}
\end{equation}
\end{example}

\begin{proposition} \cite{OSS}
\label{prop_Hbar_Cons} Let $\gf\in\{A_{n-1},C_n,D_n\}$ (so that
$\gaf\in\{A_{n-1}^{(1)},A_{2n-1}^{(2)},D_n^{(1)}\}$). For $\delta\in\Z_{>0}^m$
the embedding of $\gf$-crystals
$S_\delta^\gf:B_\delta^\gf\hookrightarrow B_{1^{|\delta|}}^\gf$
preserves coenergy, that is, $\Db^{\gf}(b)=\Db^\gf(S_\delta^\gf(b))$
for any $b\in B_\delta^\gf$.
\end{proposition}
\begin{proof} This is proved in \cite{OSS} for nonexceptional affine root systems,
including $D_n\subset D_n^{(1)}$ with the standard Dynkin labeling. The
result for $D_n^\dagger\subset D_n^{(1)\dagger}$ follows by applying the affine Dynkin automorphism
given by $i\mapsto n-i$ and the automorphism of $B_1^{D_n}$ given by $j\mapsto \overline{n+1-j}$
and $\overline{j}\mapsto n+1-j$ for $1\le j\le n$.
\end{proof}

\begin{example} Let us verify Proposition \ref{prop_Hbar_Cons} for
$\delta=(l,k)$ and the element $v=v_{l,k;a,b}^{C_n}$ of Remark
\ref{R:tensor2}. On one hand, $\Db(b)=2a+b$ by Example \ref{exa1}.
By Example \ref{ex_spli} $S^{C_n}_{l,k}(v)=1^{\otimes k+l-a-b}
\otimes 2^{\otimes b} \otimes \ba^{\otimes a}$. Computing
$\Db(S^{C_n}_{l,k}(v))$ using \eqref{E:Dhomog} and \eqref{H_simple}
we obtain the answer $(a+b)+a$, the sum of the positions (counting
from the right) of the ascents $1\otimes 2$ and $2\otimes \bar{1}$
in $S^{C_n}_{l,k}(v)$.
\end{example}

\subsection{One-dimensional sums for the affine crystals
$B_\mu^{C_n}$} For the sake of completeness we prove Lemma
\ref{L:stable} for the affine family $\gaf_n=A_{2n-1}^{(2)}$.

\begin{lemma}
\label{lem_inde} Let $m\le n$ be nonnegative integers and
$\delta\in\mathbb{N}^m$. Then every highest weight vertex $b\in
B_\delta^{C_n}$ contains only letters of the set
$\{1,\dotsc,m,\overline {m-1},\overline{m-2},\dotsc,\overline{1}\}.$
\end{lemma}
\begin{proof}
We fix $n$ and proceed by induction on $m$. When $m=1$,
$b=1^{\delta_1}$ and the property holds. Let $m\ge2$ and suppose the
above property holds for $m-1$. Let $b=b_1\otimes\dotsm\otimes
b_m\in B_\delta^{C_n}$ be a highest weight vertex for
$\delta\in\Z_{>0}^m$. By Lemma \ref{lem_plu_hp},
$b^\flat=b_1\otimes\dotsm\otimes b_{m-1}\in B_{\delta^\flat}^{C_n}$
is a highest weight vertex with
$\delta^\flat=(\delta_1,\dotsc,\delta_{n-1})$. By induction
$b^\flat$ contains only letters in
$\{1,\dotsc,m-1,\overline{m-2},\dotsc,\overline{1}\}$; in particular
$\vp_i(b^\flat)=0$ for all $i\in\{m,\dotsc,n\}$. If $b_m$
contains a letter not in the set $\{1,\dotsc,m,\overline{m-1},\overline{m-2}%
,\dotsc,\overline{1}\}$ then there exists an integer
$i\in\{m,\dotsc,n\}$ such that $\ve_i(b_m)\ge1$. By \eqref{E:tensor}
we obtain $\ve_i(b)\ge1$ which contradicts the assumption that $b$
is a highest weight vertex.
\end{proof}

\begin{corollary} \label{C:stable} Fix $m$ and $\la,\mu\in\Par_m$. Then for all $n\ge
m$ the one dimensional sum $X_{\la,\mu}^{A_{2n-1}^{(2)}}(q)$ does
not depend on $n$.
\end{corollary}
\begin{proof}
By the previous lemma the set $F_{\la,\mu}$ does not depend on $n$.
The $C_n$ plactic relations are stable for increasing $n$, so the
computations of combinatorial $R$-matrices and the local energy
functions that occur in the computation of $\Db$, are independent of
$n$. The corollary follows.
\end{proof}

\begin{remark} This also follows from the corresponding stability
result for the fermionic formula $M$ \cite{SZ} and the $X=M$ theorem
for symmetric powers \cite{SCS}, but this is massive overkill.
\end{remark}

\section{Proof of $X=K$ for $\dia=(1,1)$}
\label{sec_proof1}

Fix integers $0\le m<n$.

\subsection{Combinatorial description of $X=K$} \label{S:q=1}%
We reformulate the desired $X=K$ identity \eqref{E:X=K} for
$\dia=(1,1)$. We recall a classical result of Littlewood
\cite[appendix p. 295]{Littl} \cite[Theorem $A_1$]{KT}:
\begin{equation} \label{E:LRsumRes}
\sum_{\gamma\in \Par_m^{(1,1)}}
  c^\nu_{\la\gamma} = [V^{D_n}(\nu):V^{A_{n-1}}(\la)]\qquad\text{for
  $\nu,\la\in \Par_m$ with $n\ge m$,}
\end{equation}
which is the branching multiplicity of $V^{A_{n-1}}(\la)$ in the
restriction of $V^{D_n}(\nu)$ to $U_q(A_{n-1})$. Using
\eqref{E:LRsumRes} and Lemma \ref{L:stable}, the $X=K$ formula
\eqref{E:X=K} for $\dia=(1,1)$ may be rewritten as
\begin{equation} \label{E:X=K(1,1)}
  \Xb_{\la,\mu}^{A_{2n-1}^{(2)}}(q) = q^{\frac{|\mu|-|\la|}{2}}
  \sum_{\nu\in\Par_m} \Kb_{\nu,\mu}(q)\,
  [V^{D_n}(\nu):V^{A_{n-1}}(\la)].
\end{equation}

Let $E_{\la,\mu}$ be the set of $A_{n-1}$-highest weight vertices of
weight $\la$ in $B_\mu^{D_n}$ belonging to the image of the map
$\Psi$ in Lemma \ref{lemm_B(nu)}. Given $b\in F_\mu^{A_{n-1}}$, let
$E_{\la,\mu,b}\subset B_\mu^{D_n}$ be the set of $A_{n-1}$-highest
weight vertices of weight $\la$ in the component $B^{D_n}(b) \subset
B_\mu^{D_n}$. Then $E_{\la,\mu} = \bigcup_{b\in F_\mu^{A_{n-1}}}
E_{\la,\mu,b}$.

The identity \eqref{E:X=K(1,1)} can be rewritten as
\begin{equation} \label{E:XKcomb}
  \sum_{v\in F_{\la,\mu}} q^{\Db(v)} =
  q^{\frac{|\mu|-|\la|}{2}} \sum_{b\in F_\mu^{A_{n-1}}}
  q^{\Db^A(b)} \sum_{c\in E_{\la,\mu,b}} 1,
\end{equation}
due to the definition of the 1-d sums $\Xb_{\la,\mu}^{A_{2n-1}}(q)$
and $\Kb_{\nu,\mu}(q) = \Xb^{A_{n-1}^{(1)}}_{\nu,\mu}(q)$, and Lemma
\ref{lemm_B(nu)}.

\begin{remark} \label{R:bij} To prove Theorem \ref{T:X=K}
for $\dia=(1,1)$ it suffices to exhibit a bijection
$\theta:F_{\la,\mu}\rightarrow E_{\la,\mu}$ such that, for every
$v\in F_{\la,\mu}$ such that $\theta(v) \in E_{\la,\mu,b}$,
$\Db(v)=\Db^A(b)+\frac{|\mu|-|\la|}{2}$. It also suffices to give a
bijection $\theta:F_{\la,\mu}\rightarrow E_{\la,\mu}$ and a
statistic $\KK:B_\mu^{D_n}\rightarrow\Z$ such that
\begin{enumerate}
\item $\KK$ is constant on $D_n$-components.
\item $\Db(v)=\KK(\theta(v))+\frac{|\mu|-|\la|}{2}$.
\item $\KK(b)=\Db^A(b)$ for all $b\in F_\mu^{A_{n-1}}$.
\end{enumerate}
\end{remark}

\begin{lemma} \label{L:EF}
The sets $F_{\la,\mu}$ and $E_{\la,\mu}$ have the same
cardinality.
\end{lemma}
\begin{proof} It is equivalent to show that \eqref{E:X=K} holds for $\dia=(1,1)$ at
$q=1$. Let $s_\la^{(1,1)}$ denote the universal character of Koike
and Terada for the orthogonal groups \cite{KT} \cite{SZ}. We have
\begin{equation*}
\begin{split}
 \Xb_{\la,\mu}^{(1,1)}(1) &= [ V^{D_n}(\mu_1\omega_1)
 \otimes\dotsm\otimes V^{D_n}(\mu_m\omega_1) : V^{D_n}(\la)] \\
 &= [ s_\la^{(1,1)} ] s_{\mu_1}^{(1,1)} \dotsm s_{\mu_m}^{(1,1)} \\
 &= [ s_\la^{(1,1)} ] s_{\mu_1} \dotsm s_{\mu_m} \\
 &= [ s_\la^{(1,1)} ] \sum_{\nu\in\Par_m} K_{\nu,\mu} \,s_\nu \\
 &= \sum_{\nu\in\Par_m} \Kb_{\nu\mu}(1) \sum_{\tau\in \Par_m^{(1,1)}}
 c^\nu_{\la,\tau}
\end{split}
\end{equation*}
using the following facts: by definition, the stable 1-d sum of kind
$(1,1)$ is a graded tensor product for the affine root system
$D_n^{(1)}$ which has $D_n$ as classical subalgebra with $B_s \cong
B(s\omega_1)$; the orthogonal universal characters multiply like
irreducible $D_n$ characters in sufficiently large rank \cite{KT};
$s_r^{(1,1)}$ is the Schur function $s_r$ \cite{KT}; the
multiplicities of products of single row Schur functions are Kostka
numbers; the Kostka polynomial at $q=1$ is the Kostka number; and
the Littlewood formula for the coefficient of an orthogonal
universal character in a Schur function \cite{KT}.
\end{proof}

\subsection{$A_{n-1}$-subcrystals of $B_\mu^{C_n}$ and $B_\mu^{D_n}$}

For $s\in\Z_{>0}$ the $C_n$-crystal $B_s^{C_n}$ decomposes into
$A_{n-1}$-connected components obtained by deleting the arrows of
color $n$. This is given explicitly by
\begin{equation*}
  B_s^{C_n} \cong \bigoplus_{\substack{\alpha,\beta,\ve\ge0 \\
  \alpha+2\ve+\beta=s}} B^A(\bn^{\alpha+\ve}n^\ve 1^\beta)
\end{equation*}
where $B^A(b)$ denotes the $A_{n-1}$-connected component of $b$. For
$\alpha,\beta\ge0$ with $\alpha+\beta=s$ define
$v_{\alpha,\beta}^{C_n}=\bn^\alpha1^\beta\in B_s^{C_n}$. Define the
$A_{n-1}$-subcrystal $\Bh_s^{C_n} \subset B_s^{C_n}$ by
\begin{equation} \label{E:Chat}
  \Bh_s^{C_n} = \bigoplus_{\alpha+\beta=s} B^A(v_{\alpha,\beta}^{C_n}).
\end{equation}
$\Bh_s^{C_n}$ can also be characterized by the property that it is
the largest $A_{n-1}$-subcrystal of $B_s^{C_n}$ such that none of
its vertices contain the pair $\bn n$.

Similarly the $D_n^\dagger$-crystal $B_s^{D_n}$ decomposes as an
$A_{n-1}$-crystal as follows. For $\alpha,\beta\ge0$ such that
$\alpha+\beta=s$, let $v_{\alpha,\beta}^{D_n} =
1^\beta\bn^{\alpha}$. We have the $A_{n-1}$-decomposition
\begin{equation}
  B_s^{D_n} \cong \bigoplus_{\alpha+\beta=s}
  B^A(v_{\alpha,\beta}^{D_n}).
\end{equation}
As $A_{n-1}$ crystals, $B^A(v_{\alpha,\beta}^{C_n}) \cong
B^A(v_{\alpha,\beta}^{D_n}) \cong
B^A(\beta\omega_1+\alpha\omega_{n-1})$. Let $\theta_{\alpha,\beta}$
be the $A_{n-1}$-crystal isomorphism
\[
\theta_{\alpha,\beta}:B^A(v_{\alpha,\beta}^{C_n}) \rightarrow
B^A(v_{\alpha,\beta}^{D_n}).
\]
These maps patch together to define an $A_{n-1}$-crystal isomorphism
$\theta_s:\Bh_s^{C_n} \rightarrow B_s^{D_n}$.

\begin{remark} \label{R:theta}
\begin{enumerate}
\item In fact $B_s^{C_n} \cong
\bigoplus_{r=0}^{\lfloor\frac{s}{2}\rfloor} B_{s-2r}^{D_n}$ as
$A_{n-1}$-crystals. The element $b\in B_s^{C_n}$ is sent to the
$r$-th summand $B_{s-2r}^{D_n}$ if $r$ is the maximum of the number
of letters $n$ and $\bn$ in $b$. In particular, for $s=1$,
$\Bh_1^{C_n} \cong B_1^{C_n}$ and $\theta_1:B_1^{C_n}\rightarrow
B_1^{D_n}$ is a $A_{n-1}$-crystal isomorphism given by the identity
map.
\item The vertices of $B(v_{\alpha,\beta}^{C_n})$ (resp.
$B(v_{\alpha,\beta}^{D_n})$) do not contain any pair $(\bn,n)$
(resp. $(1,\overline{1}))$.
\item The map $\theta_{\alpha,\beta}$ may be computed applying transformations of
the form
\begin{equation}
\label{compu_teta}%
\overline{q}p\mapsto p\overline{q}\quad\text{ if $p\ne q$} \qquad
\text{ and } \qquad \overline{p}p\mapsto
(p+1)\overline{p+1}\quad\text{ if $p\ne n$.}
\end{equation}
\end{enumerate}
\end{remark}

Let $\delta\in\Z_{>0}^m$. Define the $A_{n-1}$-subcrystal
$\Bh_\delta^{C_n}\subset B_\delta^{C_n}$ by $\Bh_\delta^{C_n} =
\Bh_{\delta_1}^{C_n} \otimes\dotsm \otimes \Bh_{\delta_m}^{C_n}$ and
the map $\theta_\delta:\Bh_\delta^{C_n} \rightarrow B_\delta^{D_n}$
by $\theta_\delta = \theta_{\delta_1} \otimes\dotsm\otimes
\theta_{\delta_m}$.

\begin{lemma} \label{L:theta} The map
$\theta_\delta$ is an isomorphism of $A_{n-1}$-crystals. In
particular, for $\delta=(1^m)$, $\Bh_{1^m}^{C_n}=B_{1^m}^{C_n}$ and
$\theta_{1^m}$ is the identity, or equivalently, the restriction to
$A_{n-1}$ of the $C_n$- and $D_n$-crystal structure on the set
$B_{1^m}^{C_n}=B_{1^m}^{D_n}$ is the same.
\end{lemma}
\begin{proof} This follows immediately from Remark \ref{R:theta} (1) and
the fact that $\theta_s$ is an isomorphism of $A_{n-1}$-crystals.
\end{proof}

\begin{lemma} \label{L:thetaA} The map
$\theta_\delta$ restricts to the identity on $B_\delta^{A_{n-1}}$.
\end{lemma}
\begin{proof} It suffices to show this for one tensor factor. Let
$\delta=(s)$. The restriction of $\theta_s$ to $B_s^{A_{n-1}}$ is
the map $\theta_{0,s}:B^A(v_{0,s}^{C_n})\rightarrow
B^A(v_{0,s}^{D_n})$, which is the identity map on $B_s^{A_{n-1}}$.
\end{proof}

\subsection{The bijection $\theta:F\rightarrow E$}

\begin{proposition} \label{P:F=Estd} For $\la\in\Par_m$, $F_{\la,1^m} = E_{\la,1^m}$.
\end{proposition}
\begin{proof} By Lemma \ref{L:EF} it suffices to show that
$F_{\la,1^m}\subset E_{\la,1^m}$. Let
$b=x_1\otimes\dotsm\otimes x_m\in F_{\la,1^m}$. By restriction from
$C_n$ to $A_{n-1}$ it follows that $b\in B_{1^m}^{C_n}$ is a
$A_{n-1}$ highest weight vertex of weight $\la$. By Lemma
\ref{L:theta}, $b\in B_{1^m}^{D_n}$ is an $A_{n-1}$ highest weight
vertex of weight $\la$.

We must show that the component $B^{D_n}(b)$ is in the image of the
map $\Psi$. By Remark \ref{R:T} it suffices to show that the tableau
$P^{D_n}(x_1\dotsm x_m)$ belongs to $\TT^{D_m}$ and no contraction
happens during the corresponding insertion procedure. We proceed by
induction on $m$. When $m=1$, $b=1$ and the assertion is clear. When
$m=2$ the highest weight vertices of $B_{1^2}^{C_n}$ are
$b_1=1\otimes1,$ $b_2=1\otimes2$ and $b_3=1\otimes\overline{1}$. We
obtain
\[
P^{D_n}(11) = \young(11)\qquad P^{D_n}(12)=\young(1,2)\qquad
P^{D_n}(1\ba)=\young(1,\ba)
\]
There is no contraction during any of these insertions and the three
resulting tableaux belong to $\TT^{D_2}$. Now suppose the assertion
holds for $m-1$ with $m>2$. By Lemma \ref{lem_plu_hp} $b$ has the
form $b=b'\otimes x_m$ where $b'=x_1\otimes\dotsm\otimes x_{m-1}$ is
a highest weight vertex of $B_{1^{m-1}}^{C_n}$ and $x_m\in
B_1^{C_n}$. By Lemma \ref{lem_inde} $b'$
contains only letters of $\{1,\dotsc,m-1,\overline{m-2}%
,\dotsc,\overline{1}\}$. By the induction hypothesis,
$T'=P^{D_n}(x_1\dotsm x_{m-1})\in \TT^{D_{m-1}}$ and no contraction
occurs during the insertion procedure. By Remark \ref{R:T}(2) no
contraction occurs during the $D_n$-insertion of $x_m$ into $T'$.
Moreover we must have
$x_m\in\{1,\dotsc,m,\overline{m-1},\dotsc,\overline{1}\}$ by Lemma
\ref{lem_inde}. Thus
$x_m\overset{D_n}{\longrightarrow}T'\in\TT^{D_m}$ as desired.
\end{proof}

\begin{proposition}
\label{P:FtoE} Let $\la,\mu\in\Par_m$. Then $F_{\la,\mu}\subset
\Bh_\mu^{C_n}$ and the $A_{n-1}$-crystal isomorphism
$\theta_\mu:\widehat{B}_\mu^{C_n}\rightarrow B_\mu^{D_n}$ restricts
to a bijection $\theta:F_{\la,\mu}\rightarrow E_{\la,\mu}$.
\end{proposition}
\begin{proof} Lemma \ref{lem_inde} says that the vertices of
$F_{\la,\mu}$ do not contain any pair of letters $(\bn,n)$. By the
characterization of $\Bh_s^{C_n} \subset B_s^{C_n}$ it follows that
$F_{\la,\mu} \subset \Bh_\mu^{C_n}$.

Due to Lemma \ref{L:EF} and the injectivity of $\theta$, it suffices
to show that $\theta(F_{\la,\mu})\subset E_{\la,\mu}$. Let
$b=b_1\otimes\dotsm\otimes b_m\in F_{\la,\mu}$. Then $\theta(b)$ is
a $A_{n-1}$-highest weight vertex of weight $\la$ since $\theta$ is
an $A_{n-1}$-crystal isomorphism. By Remark \ref{R:T}(3) it suffices
to prove that the tableau $P^{D_n}(b_1\dotsm b_m)$ belongs to
$\TT^{D_m}$ and no contraction occurs during the corresponding
insertion procedure. The proof is essentially the same as in
Proposition \ref{P:F=Estd}.
\end{proof}

\begin{proposition}
\label{prop_S} For $\la\in \Par_m$ and $\delta\in\Z_{>0}^m$
the following diagram commutes:%
\begin{equation*}
\begin{CD}
F_{\la,\delta} @>{\theta_\delta}>> E_{\la,\delta} \\
@V{S_\delta^C}VV @VV{S_\delta^D}V \\
F_{\la,1^{|\delta|}} @>>{\theta_{1^{|\delta|}}}>
E_{\la,1^{|\delta|}}
\end{CD}
\end{equation*}
\end{proposition}
\begin{proof} By definition \eqref{E:Sdef} the map
$S_\delta$ can be computed using two operations: splitting the
leftmost tensor factor and swapping adjacent tensor factors using
maps of the form $B_1 \otimes B_s \rightarrow B_s \otimes B_1$ given
in Lemma \ref{ele_cry_iso}, which we shall denote by $\psi_s$.
Therefore it suffices to show that the maps $\theta$ are compatible
with these operations.

Let $b\in B_\delta^{C_n}$ be a highest weight vertex. By Lemma
\ref{lem_plu_hp} we can write $b=b_1 \otimes\dotsm\otimes b_m$ with
$b_1=1^{\delta_1}$. By direct computation we obtain the equality
$\theta_{\delta'}\circ S(b)=S\circ\theta_\delta(b)$ where $S$ and
$\delta'$ are defined in \eqref{E:S}. This is the required
compatibility of $\theta$ with splitting.

The required compatibility of $\theta$ with the swapping operation
is given by
\begin{equation} \label{E:Rtheta}
\theta_{(s,1)} \circ\psi_{1,s}^C(x\otimes
v)=\psi_{1,s}^D\circ\theta_{(1,s)}(x\otimes v)
\end{equation}
for any vertex $x\otimes v$ in $B_1^{C_n}\otimes B_s^{C_n}$ such
that there exists $\delta\in\Z_{>0}^m,$ $b=b_1\otimes\dotsm\otimes
b_m$ a highest weight vertex of $B_\delta^{C_n}$ and $1\le i\le m$
satisfying $x=b_i$ and $v=b_{i+1}$. Assume $b,x,v,i$ satisfy these
conditions. We have $v\in \Bh_s^{C_n}$ by Proposition \ref{P:FtoE}.
Let $\alpha,\beta\ge0$ be such that $v=B(v_{\alpha,\beta}^{C_n})$.
By Lemma \ref{lem_plu_hp}, as a $A_{n-1}$-crystal, $B_1^{C_n}\otimes
B(v_{\alpha,\beta}^{C_n})$ decomposes into six connected components.
Below we give the six associated highest weight vertices and their
images under $\psi_s^C$, computed using $C_n$-insertion.
\begin{align*}
b_1&=1\otimes\bn^\alpha1^\beta& \psi_l^C(b_1)&=\bn^{\alpha-1}1^{\beta+1}\otimes\bn \\
b_2&=1\otimes\overline{1}\bn^{\alpha-1}1^\beta&\psi_l^C(b_2)&=\bn^{\alpha-1}1^{\beta+1}\otimes
\overline{1} \\
b_3&=1\otimes\bn^\alpha21^{\beta-1}&\psi_l^C(b_3)&=\bn^{\alpha}1^{\beta}\otimes2 \\
b_4&=\bn\otimes\bn^\alpha1^\beta&\psi_l^C(b_4)&=\bn^{\alpha+1}1^{\beta-1}\otimes1 \\
b_5&=\bn\otimes\bnm\bn^{\alpha-1}1^\beta
&\psi_l^C(b_5)&=\bn^\alpha1^\beta\otimes\bnm \\
b_6&=\bn\otimes\bnm\bn^{\alpha-1}(n-1)1^{\beta-1}&\psi_l^C(b_6)&=\bn^{\alpha
-1}\bn n1^\beta\otimes\bn%
\end{align*}
If $x\otimes v$ belongs to the connected component of
$B_1^{C_n}\otimes B(v_{\alpha,\beta}^{C_n})$ with highest weight
vertex $b_6$, then $\sigma_i(b)$ (see \eqref{sigmai}) contains a
pair of letters $(n,\bn)$ and is a highest weight vertex of
$B_\delta^{C_n}$. This contradicts Proposition \ref{P:FtoE}. Since
$\theta,\psi_s^C$ and $\psi_s^D$ are isomorphisms of
$A_{n-1}$-crystals, it is enough to establish \eqref{E:Rtheta} for
the vertices $b_i$, $i\in\{1,\dotsc,5\}$. For
$b_1$ we obtain the commuting diagram %
\begin{equation*}
\begin{CD}
1\otimes\bn^\alpha1^\beta @>\theta>>
1\otimes1^\beta\bn^\alpha \\
@V{\psi_s^C}VV @VV{\psi_s^D}V \\
\bn^{\alpha-1}1^{\beta+1}\otimes\bn @>>\theta> 1^{\beta+1}\bn^{\alpha-1}%
\otimes\bn
\end{CD}
\end{equation*}
where $\psi_s^D$ is computed using $D_n$-insertion and $\theta$ is
obtained by \eqref{compu_teta}. For $b_1$ we only need to apply
plactic relations that preserve the letters during the $C_n$- and
$D_n$-insertion, so that \eqref{E:Rtheta} is immediate. The proof is
similar for $b_3$, $b_4$, and $b_5$. It only remains to consider
$b_2$. We must apply a plactic relation of type $R_3^D$ to compute
$\psi_s^D(b_2)$. This gives the commuting diagram
\begin{equation*}
\begin{CD}
1\otimes\overline{1}\bn^{\alpha-1}1^\beta @>\theta>>
1\otimes1^{\beta-1}2\overline{2}\overline {n}^{\alpha-1} \\
@V{\psi_s^C}VV @VV{\psi_s^D}V \\
\bn^{\alpha-1}1^{\beta+1}\otimes\overline{1} @>>\theta> 1^{\beta+1}\bn^{\alpha-1}%
\otimes\overline{1}
\end{CD}
\end{equation*}
\end{proof}

\subsection{The map $\protect\theta _{\protect\mu }$ and the coenergy
functions $\overline{D}$ and $\widetilde{D}$}

\label{S:K}

We shall establish the crucial relation
\begin{equation*}
\overline{D}(b)=\frac{|\mu |-|\lambda |}{2}+\widetilde{D}(\theta (b))
\end{equation*}
for any vertex $b\in F_{\lambda \mu }$ where $\overline{D}$ and $\widetilde{D}$
are the coenergy functions defined on $B_{\mu }^{C_{n}}$ and $B_{\mu}^{D_{n}},$ respectively.
\ We first begin with the case $\mu =(1^{m})$.
Recall that $\overline{D}\!:B_{1^{m}}^{C_{n}}\rightarrow \mathbb{Z}$
and $\widetilde{D}:B_{1^{m}}^{D_{n}}\rightarrow \mathbb{Z}$ are then the
statistics on vertices $b=x_{1}\otimes \dotsm \otimes x_{m}$ with $x_{i}\in
B_{1}^{C_{n}}=B_{1}^{D_{n}}$ defined by
\begin{align*}
\overline{D}\!(b)& =\sum_{i=1}^{m-1}(m-i)\,\overline{H}^{C}\!(x_{i}\otimes
x_{i+1}) & \qquad \widetilde{D}(b)& =\sum_{i=1}^{m-1}(m-i)\,\widetilde{H}
(x_{i}\otimes x_{i+1}) \\
\overline{H}\!(x\otimes y)& =
\begin{cases}
0 & \text{if $x\geq ^{C_{n}}y$} \\
1 & \text{if }x<^{C_{n}}y\text{ and }(x,y)\neq (1,\overline{1}) \\
2 & \text{if }(x,y)\neq (1,\overline{1}).
\end{cases}
& \qquad \widetilde{H}(x\otimes y)& =
\begin{cases}
0 & \text{if $x\geq ^{D_{n}}y$} \\
1 & \text{if $x\not\geq ^{D_{n}}y$ and $(x,y)\not=(\overline{n},n)$} \\
2 & \text{if $(x,y)=(\overline{n},n)$.}
\end{cases}
\end{align*}
In particular $\widetilde{H}(1\otimes \overline{1})=\widetilde{H}(\overline{1}\otimes 1)=1$,
$\widetilde{H}(p\otimes \overline{q})=0$ and $\widetilde{H}(%
\overline{q}\otimes p)=1$ when $p,q\in \{1,\dotsc ,n\}$ and $(p,q)\neq (1,1)$.

\begin{remark}
\label{R:K}The restriction of the map $\overline{D}$ and $\widetilde{D}$ to
$B_{1^{m}}^{A_{n-1}}$ is the coenergy function $\overline{D}^{A}$ for $
A_{n-1}^{(1)}$-affine crystals.
\end{remark}

\begin{proposition}
\label{lem_proK}

\begin{enumerate}
\item  The statistic $\overline{D}$ is constant on the $C_n$-components of $B_{1^{m}}^{C_{n}}$.

\item  The statistic $\widetilde{D}$ is constant on the $D_n^\dagger$-components of $B_{1^m}^{D_n}$.
\end{enumerate}
\end{proposition}

\begin{proof}
The two assertions hold since $\overline{D}$ and $\widetilde{D}$ are the coenergy
functions on affine crystals for types $C_n\subset A_{2n-1}^{(2)}$ and $D_n^\dagger\subset D_n^{(1)\dagger}$
respectively.
\end{proof}

For any vertex $b=x_{1}\otimes \dotsm \otimes x_{m}\in B_{1^{m}}^{C_{n}}$
let
\begin{equation*}
Z_{b}=\{i\in \{1,\dotsc ,m-1\}\mid \text{$x_{i}$ and $x_{i+1}$ are not
simultaneously barred or unbarred. }\}
\end{equation*}
Given $x\otimes y$ in $B_{1^{2}}^{C_{n}}$, we set
\begin{equation*}
\overline{h}(x\otimes y)=
\begin{cases}
0 & \text{if $x\geq ^{C_{n}}y$} \\
1 & \text{if }x<^{C_{n}}y
\end{cases}.
\end{equation*}
This means that $\overline{h}$ is the coenergy function on $B_{1^{2}}^{C_{n}}
$ regarded as an affine $C_{n}^{(1)}$-crystal. In particular we have
$\overline{h}(x\otimes y)=\overline{H}(x\otimes y)$ if $(x,y)\neq (1,\overline{1})$
and $\overline{h}(1\otimes \overline{1})=1$.

\begin{lemma}
\label{lem2} \cite{lec2} For $\lambda \in \mathcal{P}_{m}$ and any $b\in
F_{\lambda ,1^{m}}$ we have
\begin{equation}
\sum_{i\in Z_{b}}(m-i)(-1+2\,\overline{h}\!(x_{i}\otimes x_{i+1})))=\dfrac{m-|\lambda |}{2}.  \label{E:Z}
\end{equation}
\end{lemma}

\begin{proof}
We have
\begin{equation*}
-1+2\,\overline{h}(x_{i}\otimes x_{i+1})=
\begin{cases}
1 & \text{if $x_{i}<x_{i+1}$} \\
-1 & \text{otherwise.}
\end{cases}
\end{equation*}
Since $b$ is a highest weight vertex, $x_{1}=1$ by Lemma \ref{lem_plu_hp}.
We obtain
\begin{equation}
\sum_{i\in Z_{b}}(-1+2\,\overline{h}(x_{i}\otimes x_{i+1}))=
\begin{cases}
0 & \text{if $x_{m}$ is unbarred} \\
1 & \text{otherwise.}
\end{cases}
\label{s_aux}
\end{equation}
To prove the lemma we proceed by induction on $m$. When $m=1$ we have $b=1$,
$Z_{b}$ is empty, $|\la|=1$ and the lemma holds. Suppose the lemma holds for
any highest weight vertex of $B_{1^{m-1}}^{C_{n}}$. Let $b=b^{\prime}\otimes x_{m}\in B_{1^{m}}^{C_{n}}$
be a highest weight vertex of weight $%
\la$ with $b^{\prime }=x_{1}\otimes \dotsm \otimes x_{m-1}\in
B_{1^{m-1}}^{C_{n}}$, necessarily a highest weight vertex, of weight $\gamma$, say.
Let $s$ be the left hand side of \eqref{E:Z}. We have $%
Z_{b}=Z_{b^{\prime }}$ if $x_{m-1}$ and $x_{m}$ are simultaneously barred or
unbarred and $Z_{b}=Z_{b^{\prime }}\cup \{m-1\}$ otherwise. Considering the
various cases and using the induction hypothesis and \eqref{s_aux} we have
\begin{equation*}
s=\frac{m-1-|\gamma |}{2}+
\begin{cases}
0 & \text{if $x_{m}$ is unbarred} \\
1 & \text{otherwise.}
\end{cases}
\end{equation*}
When $x_{m}$ is unbarred $|\gamma |=|\la|-1$ and when $x_{m}$ is barred, $|\gamma |=|\la|+1$.
Either way we obtain $s=\frac{m-|\la|}{2}$ as desired.
\end{proof}

\begin{proposition}
\label{prop_hK} For any $b\in F_{\lambda,1^m}$ we have
\begin{equation*}
\overline{D}(b) = \widetilde{D}(b) + \frac{m-|\lambda|}{2}.
\end{equation*}
\end{proposition}

\begin{proof}
Write $b=x_{1}\otimes \dotsm \otimes x_{m}$ and let $N_{b}=\{i\in Z_{b}\mid
\text{$x_{i}=1$ and $x_{i+1}=\overline{1}$ }\}$. By Lemma \ref{lem_inde}, we
have $x_{i}\notin \{\overline{n},n\}$ for any $i=1,\ldots ,m$. This gives
\begin{equation*}
\widetilde{H}(x_{i}\otimes x_{i+1})=
\begin{cases}
\overline{H}(x_{i}\otimes x_{i+1}) & \text{if $i\not\in Z_{b}$}
 \\ 1-\overline{H}(x_{i}\otimes x_{i+1}) & \text{if $i\in Z_{b}$ and $i\not\in
N_{b}$} \\
\overline{H}(x_{i}\otimes x_{i+1})-1 & \text{if }i\in N_{b}.
\end{cases}
\end{equation*}
and using Lemma \ref{lem2} we obtain
\begin{align*}
\overline{D}(b)& =\sum_{i=1}^{m-1}(m-i)\,\overline{H}(x_{i}\otimes
x_{i+1})=\sum_{i\not\in N_{b}}(m-i)\,\overline{H}(x_{i}\otimes
x_{i+1})+\sum_{i\in N_{b}}(m-i)2 \\
& =\sum_{i\not\in Z_{b}}(m-i)\,\widetilde{H}(x_{i}\otimes
x_{i+1})+\sum_{i\in Z_{b}-N_{b}}(m-i)\,\overline{H}(x_{i}\otimes
x_{i+1})+\sum_{i\in N_{b}}(m-i)2 \\
& =\sum_{i\notin Z_{b}}(m-i)\widetilde{H}(x_{i}\otimes x_{i+1})+\sum_{i\in
Z_{b}-N_{b}}(m-i)\widetilde{H}(x_{i}\otimes x_{i+1})+\sum_{i\in N_{b}}(m-i)+
\\
& \sum_{i\in Z_{b}-N_{b}}(m-i)(-1+2\overline{H}(x_{i}\otimes
x_{i+1}))+\sum_{i\in N_{b}}(m-i) \\
& =\widetilde{D}(b)+\sum_{i\in Z_{b}-N_{b}}(m-i)(-1+2\overline{H}(x_{i}\otimes x_{i+1}))+\sum_{i\in N_{b}}(m-i) \\
& =\widetilde{D}(b)+\sum_{i\in Z_{b}}(m-i)(-1+2\overline{h}(x_{i}\otimes
x_{i+1})) \\
& =\widetilde{D}(b)+\frac{m-|\lambda |}{2}.
\end{align*}
where the last equality follows from Lemma \ref{lem2}.
\end{proof}

By Proposition \ref{prop_Hbar_Cons}, the splitting $S_{\mu }^{D}$ preserves
the coenergy function $\widetilde{D}_{\mu }$.\ Thus we have
\begin{equation*}
\widetilde{D}(b)=\widetilde{D}(S_{\mu }^{D}(b))
\end{equation*}
for any $b\in B_{\mu }^{D_{n}}$.

\begin{proposition}
\label{prop_Hbar_Kbar} Let $\lambda,\mu\in\mathcal{P}_m$ and $b\in
F_{\lambda,\mu}$.

\begin{enumerate}
\item  $\widetilde{D}$ is constant on $D_{n}^\dagger$-components.

\item  We have $\overline{D}(b)=\frac{|\mu |-|\lambda |}{2}+\widetilde{D}(\theta (b))$.

\item  Let $b\in F_{\mu }^{A_{n-1}}$. Then $\widetilde{D}(b)=\overline{D}^{A}(b)$.
\end{enumerate}
\end{proposition}

\begin{proof}
Assertion (1) follows from the fact that $\widetilde{D}$ is the coenergy
function on $B_{\mu }^{D_{n}}.$

For (2), by Proposition \ref{prop_hK} and Proposition \ref{prop_S} we have
\begin{equation*}
\Db(b)=\Db(S_{\mu }^{C}(b))=\frac{|\mu |-|\la|}{2}+\KK(S_{\mu }^{D}\circ
\theta (b))=\frac{|\mu |-|\la|}{2}+\KK(\theta (b)).
\end{equation*}
For (3), let $b\in F_{\mu }^{A_{n-1}}$. We have
\begin{equation*}
\begin{split}
\KK(b)& =\KK(S_{\mu }^{D}(b))=\KK(S_{\mu }^{D}(\theta _{\mu }(b)) \\
& =\KK(\theta _{1^{|\mu |}}(S_{\mu }^{C}(b)))=\KK(S_{\mu }^{C}(b)) \\
& =\KK(S_{\mu }^{A}(b))=\Db^{A}(S_{\mu }^{A}(b))=\Db^{A}(b)
\end{split}
\end{equation*}
by the definition of $\KK_{\mu }$ and the following facts: $\theta $ is the
identity on $B_{\mu }^{A_{n-1}}$ (Lemma \ref{L:thetaA}), the splitting maps
intertwine $\theta $ (Proposition \ref{prop_S}), $\theta _{1^{|\mu |}}$ is
the identity (Lemma \ref{L:theta}), type C and type A splitting agree on $B_{\mu }^{A_{n-1}}$
(Remark \ref{R:SAinCD}), $\KK$ is type $A$ coenergy on $B_{1^{|\mu |}}^{A_{n-1}}$ (Remark \ref{R:K} (3)),
and that splitting
preserves coenergy (Proposition \ref{prop_Hbar_Cons}).
\end{proof}

\begin{corollary}
\label{C:XK(1,1)} Theorem \ref{T:X=K} holds for $\diamondsuit=(1,1)$.
\end{corollary}

\begin{proof}
In light of Remark \ref{R:bij}, the proof is completed by Propositions \ref
{P:FtoE} and \ref{prop_Hbar_Kbar}.
\end{proof}

\section{Proof of $K=\KLS$}
\label{S:KKL} We prove Theorem \ref{T:K=KL} with a type-independent
argument.

For our realization of each of $\gf_n\in\{A_{n-1},B_n,C_n,D_n\}$,
the set of positive roots $R^+(\gf_n)$ contain a copy of the
positive roots $R^+(A_{n-1})$ for type $A_{n-1}$. The Weyl group $W$
also contains a copy of the symmetric group, the Weyl group
$S_n=W(A_{n-1})$ of type $A_{n-1}$.

We use a number of tricks. The first, presented in section
\ref{SS:altKL}, is the well-known idea \cite{G} of viewing the
Lusztig $q$-analogues as coefficients in a generating function which
is the graded character of twisted functions on the nullcone. More
precisely we use an analogous generating function for the stable KL
polynomials. This allows us to apply the second trick in section
\ref{SS:Littlewood}, which is the recognition of Littlewood's
formulae inside the above generating function. This requires a
special property of the function $L$ of Theorem \ref{T:K=KL}. Lesser
tricks include manipulations of rational $gl_n$ characters including
contragredient duality and symmetries of tensor product
multiplicities.

\subsection{Another formulation of the stable KL polynomials}
\label{SS:altKL} We give an alternative form for the definition of
the stable KL polynomials $\KLS_{\la,\mu}^{\gf_n,L}(q)$ which is
helpful for our proof of $K=\KLS$. Let $\Z[X]$ be the polynomial
ring in a set $X=(x_1,\dotsc,x_n)$ of $n$ variables, interpreted as
the group algebra of the weight lattice of $gl_n$, where
$x_i=e^{\ve_i}$ is the exponential of the $i$-th standard basis
vector in the weight lattice $\Z^n$. Let $\DW$ (see \eqref{E:DWdef}
be the set of dominant weights. Let $J=\sum_{w\in S_n}
(-1)^{\ell(w)} w$ be the antisymmetrization operator on $\Z[X]$ and
$\rho=\rho^{A_{n-1}}=(n-1,n-2,\dotsc,1,0)$ the half-sum of positive
roots.

Let $\Dem:\Z[X]\rightarrow \Z[X]^{S_n}$ be the Demazure operator for
the longest element $w_0$ of $S_n=W(A_{n-1})$, defined by the
bialternant
\begin{equation}
  \Dem(f) = J(x^\rho)^{-1} J(x^\rho f).
\end{equation}
It sends dominant monomials to irreducible rational
$gl_n$-characters
\begin{align}
  \Dem(x^\la)&=s_\la[X] &&\text{for dominant weights
  $\la\in\DW$}
\intertext{satisfies the $\rho$-shifted antisymmetry for monomials}%
  \Dem(x^\beta) &= (-1)^w \Dem(x^{w(\beta+\rho)-\rho})&\qquad&\text{for
  $\beta\in \Z^n$}
\intertext{and is linear with respect to symmetric polynomials}%
\label{E:symlinear}%
\Dem(fg)&=f \Dem(g) &&\text{for $f\in \Z[X]^{S_n}$.}
\end{align}
The character $s_\la[X]$ is a Schur polynomial up to an integral
power of $x_1x_2\dotsm x_n$. Since $\{s_\la[X]\mid \la\in\DW \}$ is
a basis of $\Z[X]^{S_n}$ we may define
\begin{equation} \label{E:KLgf}
  \sum_{\la\in\DW} s_\la[X] \, \KLS_{\la,\mu}^{\gf_n,L}(q) =
  \Dem \left(e^\mu \prod_{\alpha\in R^+(\gf_n)}
  \dfrac{1}{1-q^{L(\alpha)}e^\alpha} \right)
\end{equation}
Multiplying both sides by the Vandermonde $J(e^\rho)$ and taking
coefficients, one sees that the two definitions of stable KL
polynomials \eqref{E:stableKL} and \eqref{E:KLgf} agree for any pair
$(\la,\mu)$ of elements in $\DW$.

\begin{example} For $n=2$, type $A_{n-1}=A_1$, and $\mu=(0,0)$, the
left hand side of \eqref{E:KLgf} reads
\begin{equation*}
  s_{(0,0)} + q \,s_{(1,-1)} + q^2 s_{(2,-2)} + \dotsm
\end{equation*}
so that $\KLS_{(r,-r),(0,0)}(q)=q^r$. Here
$s_{(r,-r)}(x_1,x_2)=(x_1x_2)^{-r} s_{(2r,0)}(x_1,x_2)$ where
$s_{(2r,0)}(x_1,x_2)=\sum_{a+b=2r} x_1^a x_2^b$ is the Schur
polynomial. By \eqref{E:KLshift} and Proposition \ref{P:stableKL},
we have the Kostka-Foulkes polynomial
$K_{(2r,0),(r,r)}(q)=\KLS_{(2r,0),(r,r)}(q)=q^r$. This agrees with
\cite{LSc1}: there is one tableau of shape $(2r,0)$ and weight
$(r,r)$, the single-rowed tableau $1^r2^r$, and it has charge $r$.
\end{example}

\subsection{The Littlewood formulae} \label{SS:Littlewood}
Here is where the choice of function $L$ enters in. With the
function $L$ as in Theorem \ref{T:K=KL}, we observe that for all
$\alpha\in R^+(\gf_n)\backslash R^+(A_{n-1})$, $L(\alpha)$ is half
the length of the root $\alpha$. This gives the left hand side of
\eqref{E:sym} a homogeneity between $q$ and $X$ that allows us to
use Littlewood's formulae \cite{Littl} \cite{mac} to immediately
obtain
\begin{equation} \label{E:sym}
\prod_{\alpha\in R^+(\gf_n)\backslash R^+(A_{n-1})}
\dfrac{1}{1-q^{L(\alpha)}e^\alpha} = \sum_{\gamma\in\Par_n^\dia}
q^{|\gamma|/2} s_\gamma[X]
\end{equation}
where $\dia=\{\vn,(1),(2),(1,1)\}$ for types $A_{n-1},B_n,C_n,D_n$
respectively. We have
\begin{equation} \label{E:Little}
\begin{split}
\sum_{\la\in\DW} s_\la[X] \, \KLS_{\la,\mu}^{\gf_n,L}(q) &=
\sum_{\gamma\in\Par_n^\dia} q^{|\gamma|/2} s_\gamma[X] \Dem
\left(e^\mu \prod_{\alpha\in R^+(A_{n-1})}
\dfrac{1}{1-q^{L(\alpha)}e^\alpha} \right) \\
&= \sum_{\gamma\in\Par_n^\dia} q^{|\gamma|/2} s_\gamma[X]
\sum_{\nu\in\DW} s_\nu[X] \,\KLS_{\nu,\mu}^{A_{n-1}}(q)
\end{split}
\end{equation}
The first equality holds because \eqref{E:sym} is
$A_{n-1}$-symmetric, so by \eqref{E:symlinear} it can be pulled out
from the operand of $\Dem$ in \eqref{E:KLgf}. The second equality
holds by \eqref{E:KLgf} for type $A_{n-1}$.

\subsection{Rational $gl_n$ characters and the proof completed}
 Let $c^\la_{\gamma\nu}$ be the $gl_n$ tensor product
coefficient, the coefficient of the irreducible $gl_n$-character
$s_\la$ in the product $s_\gamma s_\nu$, for $\la,\gamma,\nu\in\DW$.
Taking the coefficient of $s_\la$ in \eqref{E:Little} we have
\begin{equation} \label{E:LR}
\begin{split}
 \KLS_{\la,\mu}^{\gf_n,L}(q) &= \sum_{\nu\in\DW} \KLS_{\nu,\mu}^{A_{n-1}}(q)
 \sum_{\gamma\in\Par_n^\dia} q^{|\gamma|/2}
 c^\la_{\gamma\nu} \\
 &= q^{(|\la|-|\mu|)/2} \sum_{\nu\in\DW} \KLS_{\nu,\mu}^{A_{n-1}}(q)
 \sum_{\gamma\in\Par_n^\dia}
 c^\la_{\gamma\nu}.
\end{split}
\end{equation}
The last equality holds because $c^\la_{\gamma\nu}=0$ unless
$|\la|=|\gamma|+|\nu|$ and $\KLS_{\nu,\mu}^{A_{n-1}}(q)=0$ unless
$|\nu|=|\mu|$. The right hand side of \eqref{E:LR} is close to that
of \eqref{E:X=K} but there are some differences; there are a number
of appearances of dominant weights, the tensor product multiplicity
looks turned around, and the power of $q$ has been shifted and
reversed.

We now exploit the fact that we are working with rational
$gl_n$-characters. All of the above works for $\mu\in\DW$. There is
an involutive automorphism on the weight lattice $\Z^n$ given by
reversal and negation: $\beta\mapsto \beta^*:=-w_0^{A_{n-1}}\beta$.
It permutes the set of positive roots $R^+(A_{n-1})$ and preserves
the set $\DW$ of dominant weights, sending $\la\in\DW$ to
$\la^*\in\DW$, which is the highest weight of the contragredient
dual of the irreducible $gl_n$-module of highest weight $\la$. With
$M$ as in \eqref{E:hat},
\begin{equation} \label{E:hatstar}
\begin{split}
\widehat{\la} &= \la^* + (M^n) \\
\widehat{\mu} &= \mu^* + (M^n).
\end{split}
\end{equation}
We have
\begin{equation} \label{E:comp}
\begin{split}
  \KLS_{\widehat{\la},\widehat{\mu}}^{\gf_n,L}(q) &= \KLS_{\la^*,\mu^*}^{\gf_n,L}(q) \\
  &= q^{(|\la^*|-|\mu^*|)/2}
  \sum_{\nu\in\DW} \KLS_{\nu,\mu^*}^{A_{n-1}}(q) \sum_{\gamma\in\Par_n^\dia}
  c^{\la^*}_{\gamma\nu} \\
  &= q^{(|\mu|-|\la|)/2}
  \sum_{\nu\in\DW} \KLS_{\nu^*,\mu^*}^{A_{n-1}}(q) \sum_{\gamma\in\Par_n^\dia}
  c^{\la^*}_{\gamma\nu^*} \\
  &= q^{(|\mu|-|\la|)/2} \sum_{\nu\in\DW} \KLS_{\nu,\mu}^{A_{n-1}}(q) \sum_{\gamma\in\Par_n^\dia}
  c^{\nu}_{\la\gamma}.
\end{split}
\end{equation}
The first equality is given by \eqref{E:hatstar} and the invariance
property \eqref{E:KLshift}. The second holds by applying
\eqref{E:LR} for the pair $(\la^*,\mu^*)$. The third holds since
$\nu\mapsto\nu^*$ is a bijection from $\DW$ to itself. The last
equality holds by the symmetry of rational $gl_n$ tensor product
multiplicities
\begin{equation*}
c^{\la^*}_{\gamma\nu^*} = \inner{s_\la^*}{s_\gamma s_\nu^*} =
\inner{s_\nu}{s_\la s_\gamma}= c^\nu_{\la\gamma}
\end{equation*}
and the duality
\begin{equation} \label{E:KLAdual}
\KLS_{\la^*,\mu^*}^{A_{n-1}}(q) = \KLS_{\la,\mu}^{A_{n-1}}(q)
\end{equation}
which follows, for example, from applying the automorphism
$e^\beta\mapsto e^{\beta^*}$ to the definition \eqref{E:stableKL}
for $\gf_n=A_{n-1}$.

We now restrict ourselves to the original situation where $\la$ and
$\mu$ are assumed to be partitions with at most $m$ parts where
$m\le n$. We may assume that $\nu$ is a partition, because $\gamma$
and $\la$ are partitions, and in this case, $c^\nu_{\gamma\la}=0$ if
$\nu\in\DW$ is not a partition. We may assume that $\nu\ge\mu$, for
otherwise $\KLS_{\nu,\mu}^{A_{n-1}}(q)=0$; see the proof of
Proposition \ref{P:stableKL}. Since $\mu$ is a partition with at
most $m$ parts, it follows that $\nu$ is as well. The tensor product
$c^\nu_{\gamma\la}$ is now indexed by three partitions. We may
assume that $\gamma$ has at most $m$ nonzero parts, since
$c^\nu_{\gamma\la}=0$ unless $\gamma\subset\nu$.

Finally, since $\nu$ and $\mu$ are partitions we have
$\KLS_{\nu,\mu}^{A_{n-1}}(q)=K_{\nu,\mu}(q)$ by Proposition
\ref{P:stableKL}. Therefore \eqref{E:comp} becomes
\begin{equation*}
  \KLS_{\widehat{\la},\widehat{\mu}}^{\gf_n,L}(q) =
  q^{(|\mu|-|\la|)/2} \sum_{\nu\in\Par_m}  K_{\nu\mu}(q)
  \sum_{\gamma\in\Par_m^\dia} c^\nu_{\gamma\la}.
\end{equation*}
Theorem \ref{T:K=KL} follows by replacing $q$ by $q^{-1}$,
multiplying by $q^{||\mu||+|\mu|-|\la|}$ and using the definitions
\eqref{E:XKnull} and \eqref{E:Kdef} of the cocharge Kostka-Foulkes
polynomial $\Kb_{\nu\mu}(q)$ and the $K$-polynomial
$\Kb_{\la,\mu}^\dia(q)$.

\section{Combinatorial question}
Since KL polynomials are $q$-analogues of weight multiplicities, in
light of the $X=\KLS$ theorem (Corollary \ref{C:X=KL}) there is a
grade-preserving bijection from the set $F_{\la,\mu}$ that indexes
the one-dimensional sum $X$, to the elements of weight $\mu$ in the
crystal graph of the irreducible $U_q(\gf_n)$-module of highest
weight $\la$, with some grading on the latter objects. It would
interesting to find a natural statistic on the above weight vectors
and a natural bijection of this kind.

\end{document}